\documentclass[12pt,titlepage]{article}
\usepackage{graphicx}
\usepackage{latexsym}
\usepackage{amssymb}

\title{Periodic discrete conformal maps.}
\author{U Hertrich-Jeromin
\thanks{Partially supported by the Alexander von Humboldt Stiftung as well as by
NSF grant DMS-9626804.}\\
I McIntosh\thanks{Partially supported by NSF grants DMS-9626804 and DMS-
9705479.}\\ P Norman\thanks{Partially supported by NSF grant DMS-9626804.}\\
F Pedit\thanks{ Partially supported by NSF grants DMS-9626804 and DMS-9705479.}\\
\\
Center for Geometry, Analysis, Numerics and Graphics\\
Department of Mathematics\\
University of Massachusetts\\
Amherst, MA 01003, USA}
\date{December 31st, 1998}

\renewcommand{\P}{{\bf P}}
\newcommand{\Z}{{\bf Z}}
\newcommand{\R}{{\bf R}}
\newcommand{\C}{{\bf C}}
\newcommand{\caL}{\mathcal{L}}
\newcommand{\caE}{\mathcal{E}}
\newcommand{\caF}{\mathcal{F}}
\newcommand{\caK}{\mathcal{K}}
\newcommand{\caA}{\mathcal{A}}

\newcommand{\caO}{\mathcal{O}}
\newcommand{\caZ}{\mathcal{Z}}
\newcommand{\caG}{\mathcal{G}}
\newcommand{\Tr}{{\rm Tr}}
\newcommand{\Hom}{{\rm Hom}}
\newcommand{\im}{{\rm im}}
\newtheorem{thm}{Theorem}
\newtheorem{prop}{Proposition}
\newtheorem{lem}{Lemma}
\newtheorem{cor}{Corollary}

\begin{document}
\maketitle


\section{Introduction.}

Recently there has been much interest in the theory of discrete surfaces in
3-space and its connection with the discretization of soliton equations (see
e.g.\ \cite{BobP2,NijC} and references therein). In this article we study a
discrete geometry which is the simplest example for both theories. 
Following \cite{Bob1,BobP1} we will define a
discrete conformal map (DCM) to be a map $z:\Z^2\to\P^1$ with the property that
the cross-ratio of each fundamental quadrilateral is the same.
Specifically, for four points $a,b,c,d$ on $\P^1$ define their cross-ratio to be
\[
[a:b:c:d] = \frac{(a-b)(c-d)}{(b-c)(d-a)}.
\]
Then $z:\Z^2\to\P^1$ is discrete conformal when
\begin{equation}
\label{eq:dcmintro}
[z_{k,m+1}:z_{k,m}:z_{k+1,m}:z_{k+1,m+1}] = q
\end{equation}
for some constant $q\neq 0,1,\infty$ for all $k,m$.

\begin{figure}[ht]
\centering
\includegraphics[scale=0.7]{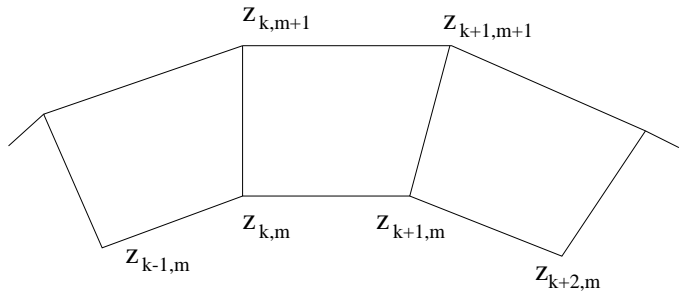}
\caption{The points about $z_{k,m}$ with neighbours joined by edges.}
\end{figure}

\noindent
The motivation for this definition is
that if $z:\R^2\to\C^1$ is smooth then it is (weakly) conformal precisely when 
\[
\lim_{\epsilon\to 0}
[z(x,y+\epsilon):z(x,y):z(x+\epsilon,y):z(x+\epsilon,y+\epsilon)] = -1,
\]
i.e.\ when $z_y^2/z_x^2 = -1$. 
Moreover, Nijhoff \& Capel \cite{NijC} have shown that one can
think of the equations (\ref{eq:dcmintro}) as being a discretization of the Schwarzian KdV (SKdV)
equations, hence
this geometry should be `integrable' in some appropriate sense. 
From another perspective, Bobenko \cite{Bob1} has shown that all the circle
patterns of Schramm \cite{Sch} correspond to DCM's with cross-ratio $-1$. 
However, as we shall see, we achieve greater insight by allowing the
cross-ratio to be any complex value. 

Our main aim is to show that all periodic DCM's (i.e.\ $z_{k+n,m}=z_{k,m}$ for
some $n$ and all $k,m$) can be constructed using methods which are straight from integrable
systems theory, viz, by relating each such map to a linear flow on the Jacobi variety 
of a compact Riemann surface $\Sigma$ (or more generally, an algebraic curve). Recall that this is 
the moduli space of degree zero holomorphic line bundles over the
curve: it is a complex manifold with the structure of an abelian group and its
dimension equals the genus of $\Sigma$. In
our case the flow is discrete so by `linear' we mean the flow is
a map of $\Z^2$ into this Jacobian which is essentially a homomorphism (generically the
map is a {\em zigzag} i.e.\ a homomorphism on a subgroup of index two, but
each of these is just a deformation of a homomorphism). 
We show that every periodic DCM is determined, uniquely up to M\" obius equivalence, by its
{\em spectral data}, which consists of: a compact hyperelliptic Riemann surface $\Sigma$ (which may be
singular) equipped with three marked points $O,S,Q$; a degree two
rational function $\lambda$ on $\Sigma$ for which $\lambda(O)=0$, $\lambda(S)=1$, $\lambda(Q)=q$;
and, a degree $g+1$ line bundle $\caL$ over $\Sigma$ satisfying a non-speciality condition, where $g$ is
the genus of $\Sigma$. 

The spectral data arises by considering the DCM as the `conformal flow' of a periodic discrete curve
i.e.\ each of the discrete periodic curves $\Gamma_{k,m} = (z_{k,m},z_{k+1,m},\ldots,z_{k+n-1,m})$ is
considered to be the evolution of the initial curve $\Gamma_{0,0}$ according to the cross-ratio
condition. Given $\Gamma_{0,0}$ and a
cross-ratio $q$ one asks the question: what condition must a point $z\in\P^1$ satisfy for it to be a
neighbour of $z_{0,0}$ in this flow? This is a question about the fixed points of 
a composite of M\" obius transformations as we go around $\Gamma_{0,0}$. We call this composite the
holonomy $H_{0,0}$ of the closed curve $\Gamma_{0,0}$. By treating the cross-ratio as a parameter 
(which we re-label $\lambda$) the holonomy becomes a rational function of $\lambda$ with values 
in $\P GL_2$. The fixed points of $H_{0,0}$ are the eigenlines of its matrix representation: these vary with
$\lambda$. The characteristic polynomial of this matrix 
determines $\Sigma$
while $\caL$ is the {\it dual} of its bundle of eigenlines. 
As a result, the $\P^1$ in which the discrete map takes values gets identified with 
the projective space
$\P\Gamma(\caL)^*$ of hyperplanes (i.e.\ dual lines) in $\Gamma(\caL)$, the space of 
globally holomorphic sections of $\caL$. 

When we do the same construction for
$\Gamma_{k,m}$ we obtain another holonomy matrix, $H_{k,m}$, with its spectral curve and line bundle
$\caL_{k,m}$. But $H_{k,m}$ is conjugate to $H_{0,0}$ by a matrix
of rational functions of $\lambda$, so the spectral curves are isomorphic. Moreover, since the conjugacy
maps eigenlines to eigenlines we obtain a rational 
section of the degree zero line bundle $\Hom(\caL_{k,m},\caL_{0,0})\simeq \caL_{0,0}\otimes\caL_{k,m}^{-1}$. 
We can explicitly compute the divisor $D_{k,m}$ of poles and zeroes of this section.
In the simplest case, 
where $\lambda=0$ is a branch point, $D_{k,m}$ is the divisor $k(S-O)+m(Q-O)$ whence the conformal flow
`linearises' on the Jacobian of $\Sigma$. 

However, the periodicity condition requires a little more:
$n(S-O)$ must be the divisor of a rational function on the singularisation $\Sigma^\prime$
of $\Sigma$ obtained by
identifying the two points over $\lambda=\infty$. This suggests that we should think of
the linearised flow as taking place on the generalised Jacobian $J^\prime$ 
for this singular curve. This leads us to
a fairly elegant formula for periodic DCM's involving the $\theta$-function for $\Sigma$ pulled back 
to $J^\prime$. This is analogous to the formula found in \cite{BobP3} for discrete surfaces of negative
Gaussian curvature.

This much is contained in sections 2 and 3. Section 2 treats the holonomy matrix for a
discrete curve and derives the spectral data. For simplicity we assume that $\Sigma$ is a non-singular
curve (we show in the Appendix that this is the generic case). Section 3 applies this to the construction
of periodic DCM's and proves that the spectral data $(\Sigma,\lambda,\caL,O,S,Q)$ characterizes the DCM
uniquely in its M\" obius equivalence class. We give the explicit formula for $z_{k,m}$ in terms of
the $\theta$-function and show that these maps will have singularities whenever the flow passes through
(a certain translate of) the $\theta$-divisor. These singularities manifest as the collapse of all four
neighbours of a point $z_{k,m}$ onto that point: in this case the cross-ratio breaks down in the adjacent
quadrilaterals. We also give a geometric interpretation for the $\theta$-function
formula which supports the view that the Schwarzian KdV equations are (one) continuum limit of the equations
(\ref{eq:dcmintro}). Geometrically this limit is very easy to describe. Let $\caA^\prime:\Sigma^\prime\to
J^\prime$ be the Abel map for $\Sigma^\prime$, then the SKdV limit arises as the secant $\vec{OS}$ 
(on $\caA^\prime(\Sigma^\prime)$)
tends to the tangent line at $O$ while $\vec{OQ}$ tends to the third derivative
$\partial^3\caA^\prime/\partial\zeta^3$ at $O$ (where $\zeta$ is a local parameter about $O$).
Finally, we compute some examples in the case where $\Sigma$ is a rational nodal curve. These
we interpret as the soliton solutions for this theory: recall that the soliton solutions of the KdV
equation have rational nodal spectral curves. Indeed, computer investigations show that the multi-soliton
solutions behave like superposed 1-solitons (in this geometry 1-solitons can be distinguished by their
rotational symmetries). 

Section 4 addresses a different aspect of the integrable nature of this geometry: the existence of the
Lax pair found by Nijhoff \& Capel \cite{NijC}. This, together with the fact that DCM's naturally occur in
1-parameter families indexed by $\lambda$, suggests that we can apply the loop group dressing action
theory known for KdV \cite{SegW,Wil} and SKdV \cite{McI}. This is achieved without much difficulty and 
we show that the appropriate dressing action preserves the cross-ratio of the map. We
concentrate on describing the corresponding `dressing orbit of the vacuum solution'. 
The term `vacuum
solution' comes from soliton theory and means the most elementary solution. In our case there is a vacuum
solution for every cross-ratio: each is the tiling 
of the plane by a given parallelogram. We show that each dressing orbit of a vacuum solution
is infinite dimensional
and their union contains every periodic DCM which is not too great a perturbation of its
continuum limit (by comparison, one knows from \cite{SegW} that every solution of KdV arising from a spectral 
curve lies in the dressing
orbit of the vacuum). The dressing construction also gives many DCM's which are not periodic 
and we compute a simple example: the discrete cubic. 

Finally, we prove a result about the connection between Darboux (B\" ack\-lund) transforms of the KdV
equation and DCM's. In \cite{Bob2} Bobenko argues that whenever a smooth integrable geometry possesses a
discrete infinite group of B\" acklund transforms this family can be taken to be an integrable
discretization of the smooth geometry. We show that each DCM in the dressing orbit of the vacuum
solution corresponds to a $\Z^2$ family of Darboux transforms of a solution of the KdV
equation. In the case of finite type (i.e.\ those DCM's possessing a spectral curve) these Darboux 
transforms are precisely the ones described in
\cite{EhlK} as those preserving the spectral curve of the KdV solution. 

\smallskip\noindent
{\bf Acknowledgements.} We thank Nick Schmitt for designing a computing package which allowed us to
investigate the moduli space of one and two solitons. Figures 5 and 6 come from this package. 

\section{Discrete Curves in $\P^1$.}

\subsection{Preliminaries.}

Let $\Gamma=(z_0,\ldots,z_{n-1})$ be distinct ordered points in $\P^1$ (and to avoid trivialities
we take $n\geq 4$). We will call
$\Gamma$ a periodic discrete curve with base point $z_0$. 
To each edge $(z_k,z_{k+1})$ (with $z_n=z_0$) we
associate a rational map of the $\lambda$-sphere $\P^1_\lambda$ (i.e.\ $\P^1$ with an
affine coordinate $\lambda$) into the M\" obius group,
\begin{equation}
\label{eq:T}
\P^1_\lambda\rightarrow \P GL_2\quad ,\quad \lambda\mapsto T_k^\lambda,
\end{equation}
by requiring that, for all $z\in\P^1\backslash\{z_k,z_{k+1}\}$, 
the cross-ratio condition
\[
[z:z_k:z_{k+1}:T_k^\lambda(z)]=\lambda
\]
is satisfied. A simple calculation shows that
$T_k^\lambda$ is invertible except at $\lambda=0,1$ and that we can represent
it in $\mathfrak{gl}_2$ by 
\begin{equation}
\label{eq:Tmatrix}
T_k^\lambda = I - \lambda^{-1}A_k
\end{equation}
where $A_k$ is the projection matrix with kernel $z_k$ and image $z_{k+1}$ (thinking of
these as lines in $\C^2$). 

\begin{figure}[ht]
\centering
\includegraphics[scale=0.6]{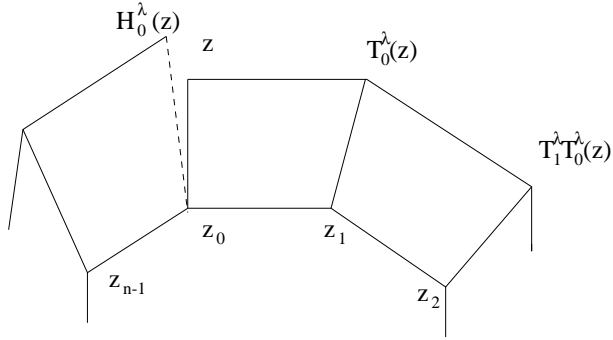}
\caption{$T_0^\lambda(z)$ and $H_0^\lambda(z)$.}
\end{figure}

\medskip\noindent
We now introduce the holonomy of $\Gamma$ at the base point $z_0$ (see Figure 2):
\begin{equation}
\label{eq:holonomy}
H_0^\lambda = T_{n-1}^\lambda\circ T_{n-2}^\lambda\circ\ldots\circ T_0^\lambda .
\end{equation}
We will use this notation for both the map into $\P GL_2$ and the matrix
representation corresponding to (\ref{eq:Tmatrix}).
\begin{lem}
\label{lem:H}
Let $d$ denote the degree of $H_0^\lambda$ (in $\lambda^{-1}$). Then $d\leq n/2$ for $n$ even and 
$d\leq (n+1)/2$ for $n$ odd.
\end{lem}
{\bf Proof.} Since $A_{k+1} A_k = 0$ the highest order term that can possibly appear in the
expansion $H_0^\lambda= I + h_1\lambda^{-1} +\ldots +h_d\lambda^{-d}$ is given by
\[
h_d = \left\{\begin{array}{cc} 
\lambda^{-n/2}(A_{n-2} B_1 A_0 + 
A_{n-1} B_2 A_1 + A_{n-1} B_3 A_0) & \hbox{\rm for $n$ even;}\\
\lambda^{-(n+1)/2}(A_{n-1} A_{n-3}\ldots A_0) & \hbox{\rm for $n$ odd;} 
\end{array}
\right.
\]
where the $B_j$ are some composites of the $A_k$.$\Box$

\subsection{Spectral Data.}

The spectral data will be computed from the trace free part of the holonomy matrix. Set
$p(\lambda)=\Tr (H_0^\lambda)/2$; then one readily sees that $p(\lambda)$ is a polynomial in
$\lambda^{-1}$ of degree $d$ with $p(\lambda)=1 + p_1\lambda^{-1} +\ldots +p_d\lambda^{-d}$. 

\smallskip\noindent
{\bf Remark.} For $n$ odd we observe from the formula above that $z_0$ is both the 
kernel and image of $h_d$, whence $h_d^2=0$ i.e.\ it is nilpotent. 
Therefore $p_d=0$ for $n$ odd.

\smallskip\noindent
Define 
\begin{equation}
\label{eq:M}
M_0^\lambda = \lambda(H_0^\lambda - p(\lambda)I)
\end{equation}
This is a trace-free matrix polynomial in $\lambda^{-1}$ of degree $\leq d-1$. 
Set $m(\lambda) = \det(M_0^\lambda)$: this is a polynomial in $\lambda^{-1}$ 
with leading order term 
$\lambda^{2-2d}\det(h_d-p_dI)$, so by the previous remark we see that
$\deg(m(\lambda))$ is at most $2d-2$ for $n$ even and $2d-3$ for $n$ odd. By the lemma
this means that for any $n$ this degree is at most $n-2$. From now on we will make a genericity
assumption: we will assume that $\deg(m(\lambda))=n-2$, that $m(\infty)\neq 0$
and that $m(\lambda)$ has distinct roots. 
Since the map $(z_0,\ldots,z_{n-1})\mapsto \det(M_0^\lambda)$ is algebraic it
is clear that the generic discrete curves occupy a Zariski open subset of
$\{(z_0,\ldots,z_{n-1}):z_i\neq z_j\}$. We will show in the appendix that it is not empty, so such
discrete curves exist and are indeed generic.
With these assumptions we have $d=n/2$ for $n$ even, $d=(n+1)/2$ for $n$ odd.

Define the spectral curve to be the isomorphism class $\Sigma$ of the curve
\[
\Sigma_0 = \{(\lambda, [v])\in \P^1\times\P^1: M_0^\lambda [v] = [v]\}.
\]
In this notation $[v]$ denotes the line through $v\in\C^2$.
\begin{prop}
\label{pp:curve}
This construction makes
$\Sigma$ a complete non-singular hyperelliptic curve of genus $g=d-2$ 
(equal to $(n-4)/2$ for $n$ even, $(n-3)/2$
for $n$ odd). This curve comes equipped with a
rational function $\lambda$ of degree two and a degree $g+1$ map $f:\Sigma\to\P^1$. The
function $\lambda$ is unbranched at $1$ and $\infty$ but when $n$ is
odd it is branched at $0$. 
\end{prop}
{\bf Proof.} By the genericity assumption $\Sigma$ is modelled by the non-singular 
completion of the affine curve with equation
$\det(\mu I-M_0^\lambda)=\mu^2+m(\lambda)=0$. This is clearly a 
hyperelliptic curve with hyperelliptic cover
$\lambda:\Sigma\to\P^1$.  Since we have assumed $m(\infty)\neq 0$ there is no
branch point at $\infty$. Since $\deg(m(\lambda)) = n-2$ this cover is branched at $\lambda=0$
(i.e.\ $\lambda^{-1} = \infty$) when $n$ is odd and we read off the genus from $\deg m(\lambda) =
2g+2$ for $n$ even and $\deg m(\lambda) = 2g+1$ for $n$ odd, giving $g=d-2$. To show $\lambda$ is
unbranched at $1$ it suffices to observe that $z_1,z_{n-1}$ are distinct eigenlines of $M_0^1$. 
To see this, simply 
note that $H_0^1 = (I-A_{n-1})\circ\ldots\circ (I-A_0)$ and $z_1=\ker (I-A_0)$ while
$z_{n-1} = \im (I-A_{n-1})$. 

In $\Sigma\times\C^2$ we have the kernel line bundle
of $\mu I - M_0^\lambda$, which is clearly holomorphic.
Its projectivisation is a holomorphic map $\Sigma\to\Sigma\times\P^1$ and we obtain
$f:\Sigma\to\P^1$ by composing this with projection on the second factor. Clearly $\Sigma_0$ is
the image of $\lambda\times f:\Sigma\to\P^1\times\P^1$. We have to show that this is an
embedding. Certainly it is injective, since the eigenlines of
$M_0^\lambda$ are distinct away from branch points of $\lambda$. Moreover
it is an embedding, for when $d\lambda = 0$ we are at ramification points,
which lie over the roots of $m(\lambda)$ (and the point over $\lambda=0$ when $n$ is odd). By the
genericity assumption at these points $M_0^\lambda$ is tranverse to the determinant conic,
whence its eigenlines have distinct tangents i.e.\ $df\neq 0$. Finally, since the image curve has
genus $g$ in the quadric surface $\P^1\times\P^1$ it must be of type $(2,g+1)$ (see e.g.\ \cite{Har}) 
i.e.\ the degree of $f$ is $g+1$. $\Box$ 

Let us define two triples $(\Sigma,\lambda,f)$ and $(\Sigma^\prime,\lambda^\prime,f^\prime)$, of data of the type in the
previous proposition, to be isomorphic if there is an
isomorphism $\Sigma\simeq\Sigma^\prime$ 
which {\em identifies} $\lambda$ with $\lambda^\prime$
and equates $f$ with a M\" obius transform of $f^\prime$. 
Then we have the following lemma.	
\begin{lem}
Any M\" obius transformation of $\Gamma$ leaves the isomorphism class of $(\Sigma,\lambda,f)$
fixed. 
\end{lem}
{\bf Proof.} Let $\Gamma^\prime= (gz_0,\ldots,gz_{n-1})$ for some M\" obius transformation $g$.
Clearly all the maps/matrices $T^\lambda,H^\lambda,M^\lambda$ are conjugated by $g$ and so we have
\[
\Sigma_0^\prime = \{(\lambda, [v])\in \P^1\times\P^1: gM_0^\lambda g^{-1}[v] = [v]\}
\]
whence the map $(\lambda,[v])\mapsto (\lambda, g[v])$ is an isomorphism between $\Sigma_0$ and
$\Sigma^\prime_0$ 
which identifies $\lambda$ with $\lambda^\prime$
and equates $f$ with a M\" obius transform of $f^\prime$.$\Box$ 

\medskip\noindent
Henceforth we will use $(\Sigma,\lambda,f)$ to denote this isomorphism class. The
{\em spectral data}  for the based curve $\Gamma$ is the quintuple $(\Sigma,\lambda,f,O,S)$ where
$O,S$ are points on $\Sigma$ such that
$(\lambda,f)(O)=(0,z_0)$ and $(\lambda,f)(S)=(1,z_1)$. That such points exist follows from:
\begin{lem}
\label{lem:zk}
The points $(0,z_0)$ and $(1,z_1)$ both lie on $\Sigma_0$.
\end{lem}
{\bf Proof.} We have already shown that $(1,z_1)$ lies on $\Sigma_0$ in the proof of the previous
proposition. Now, given the calculations in the proof of lemma \ref{lem:H} we want to show that $z_0$
is an eigenline of $h_d-\Tr(h_d)/2$. When $n$ is odd, $h_d$ is nilpotent with 
kernel $z_0$, so $(0,z_0)$ is the
ramification point over $\lambda=0$. When $n$ is even we know $h_d$ is a sum of three matrices:
the first has kernel $z_0$ and image $z_{n-2}$, the second has kernel $z_1$ and image $z_0$ while the
third has kernel $z_0$ and image $z_{n-1}$. 
It follows that $h_d$ maps the line $z_0$ to itself i.e.\ it is an eigenline.$\Box$

It is much more useful to take, in place of the M\" obius class of $f$, the line bundle
$\caL = f^*\caO_\P(1)$ i.e. the pullback of the hyperplane line
bundle over $\P^1$. Indeed, 
by definition $\caL$ is the dual to the bundle of eigenlines over $\Sigma$ 
and therefore contains exactly the information we require. It clearly has degree $g+1$ since $f$ does.
As a result of the next lemma we may recover $f$ up to isomorphism
from $\caL$ as the map  $\Sigma\to\P\Gamma(\caL)^*$ in which
$P$ is mapped to the hyperplane $\Gamma(\caL(-P))$ of all sections vanishing at $P$. 
\begin{lem}
\label{lem:embedding}
$\Gamma(\caL(-P-\tilde P))=0$ for each point $P\in\Sigma$ (where $\tilde P$
denotes the hyperelliptic involute of $P$). Thus $\dim\Gamma(\caL)=2$ and the map
$\Sigma\to\P\Gamma(\caL)^*$; $P\mapsto \Gamma(\caL(-P))$ separates points in hyperelliptic involution. 
\end{lem}
{\bf Proof.} Let $V=\P^1\times\C^2$, we will show that this is isomorphic to the direct image
$\lambda_*\caL$. It follows that
$\Gamma(\caL)=\Gamma(V)$ and a global section of $\caL$ vanishes at $P+\tilde P$ precisely when the
corresponding section of $V$ vanishes at $\lambda(P)$. Since all global sections of $V$ are constant this
will prove the lemma. Observe that the sheaf of local sections of the dual, $V^*\simeq V$, is an 
$\caO_\Sigma$-module: for
any local section $\sigma$ of $V^*$ and locally regular function $r(\lambda,\mu)$ on $\Sigma$ we define
$r(\lambda,\mu)\sigma = \sigma\circ r(\lambda I,M_0^\lambda)$. Now let $\caE\subset\Sigma\times\C^2$ denote
the eigenline bundle whose dual is $\caL$, then the natural pairing gives rise to an injective
$\caO_\Sigma$-module homomorphism of $V^*$ into $\hbox{Hom}(\caE,\caO_\Sigma)$. Therefore as an
$\caO_\Sigma$-module $V^*\simeq\caL(-D)$ for some positive divisor $D$ of degree $d$. But as an
$\caO_\P$-module $V^*\simeq\lambda_*\caL(-D)$ and $V^*$ has Euler characteristic $\chi(V^*)=2$, which gives
$\chi(\caL(-D))=2$. But in fact $D$ must be the trivial divisor, since
by Riemann-Roch $\chi(\caL(-D))=(g+1-d)+1-g$ so $d=0$, whence
$V^*\simeq\lambda_*\caL$.$\Box$

\subsection{Change of base point.}

Given the spectral data we wish ultimately to recover the discrete curve $\Gamma$. We have seen that the
point $O$ corresponds to the base point $z_0$ of $\Gamma$ via the line bundle $\caL$. 
Here we will show that the change of base point corresponds to moving only the line bundle $\caL$, not the
other spectral data.
We examine what happens when the based curve $\Gamma_0=\Gamma$
is subjected to a cyclic permutation to give $\Gamma_k=(z_k,z_{k+1},\ldots,z_{k-1})$. For
$\Gamma_k$ we have the corresponding holonomy $H_k^\lambda$ with base point $z_k$ and its trace
free part $M_k^\lambda$. Clearly we have the relationship
\begin{equation}
\label{eq:Mholonomy}
M_{k+1}^\lambda = T_k\circ M_k^\lambda\circ T_k^{-1}.
\end{equation}
Let the spectral data for $\Gamma_k$ be $(\Sigma_k,\lambda_k,\caL_k,O_k,S_k)$. In particular,
$f_k(O_k)=z_k$ and $f_k(S_k)=z_{k+1}$. Recall that when $n$ is odd $O_k$ is a ramification point of 
$\lambda$ and therefore a fixed point of the hyperelliptic involution $P\mapsto\tilde P$. 
\begin{prop}
\label{pp:flow}
For each $k$ we have an isomorphism $(\Sigma_k,\lambda_k)\simeq (\Sigma,\lambda)$ such that
$S_k$ is mapped to $S$ but $O_k$ maps to $O$ for $k$ even and $\tilde O$ for $k$ odd. Further
\[
\caL_{k+1}\otimes\caL_k^{-1}\simeq \caO_\Sigma(\tilde O_k-S). 
\]
\end{prop}
{\bf Proof.} We will construct the isomorphisms $(\Sigma_k,\lambda_k)\simeq
(\Sigma_{k+1},\lambda_{k+1})$ and then deduce the result from these. Fix $k$ and consider the map
\begin{equation}
\label{eq:Sigma}
\begin{array}{lcr} \Sigma_k & \rightarrow &\Sigma_{k+1}\\
                (\lambda,[v]) & \mapsto & (\lambda, [T_k v]) \end{array}
\end{equation}
which we deduce from (\ref{eq:Mholonomy}). Since $T_k$ is invertible except at $\lambda=0,1$ this map
is certainly biholomorphic off the points over $\lambda=0,1$ and equates $\lambda_k$ with
$\lambda_{k+1}$. Now we consider $T_k^\lambda$ about $\lambda=0,1$, where it has at most simple
zeroes or poles.

Set $\eta = \lambda^{-1} -1$: this is a local coordinate about both points $S_k,\tilde S_k$ over
$\lambda =1$. Let $v_\eta = v_0 +\eta v_1 +\ldots$ be the expansion for a locally holomorphic family
of eigenvectors for $M_k^\lambda$ about $\eta=0$. Then from $T_k^\lambda = (I-A_k)-\eta A_k$ we see
that
\[
T_k^\lambda v_\eta = (I-A_k)v_0 +\eta[(I-A_k)v_1 - A_kv_0] + O(\eta^2).
\]
If $[v_0] = z_{k+1} = \im A_k$ this has a simple zero,
otherwise it has no zero. Since we may rescale
$T_k$ without changing the map (\ref{eq:Sigma}) we see that replacing 
$T_k^\lambda$ by $\eta^{-1}T_k^\lambda$ exhibits (\ref{eq:Sigma}) as a biholomorphic map 
about $S_k=(1,z_{k+1})$. Further, to see that $S_k$ is mapped to $S_{k+1}=(1,z_{k+2})$ 
it is enough to see that
$\tilde S_k$ is not mapped to it. But $\im(I-A_k) = z_k$ so from the expression above
$\tilde S_k$ is mapped to $(1,z_k)$.

To perform a similar calculation about $\lambda=0$ we have to consider the two cases:	
$\lambda=0$ is a branch or is not a branch. In the latter case $\lambda$ is a local parameter about
each point $O_k,\tilde O_k$. Any locally holomorphic family of eigenvectors for $M_k^\lambda$ has
expansion $v_\lambda = v_0 + \lambda v_1 +\ldots$ about $\lambda=0$, whence
\[
T_k^\lambda v_\lambda = (I-\lambda^{-1}A_k)v_\lambda= -\lambda^{-1}A_kv_0 + (v_0-A_kv_1) + O(\lambda).
\]
This has a simple pole unless $[v_0] = z_k=\ker A_k$ i.e.\ a simple pole only at $\tilde O_k$. By
replacing $T_k$ with $\lambda T_k$ about $\tilde O_k$ we see that (\ref{eq:Sigma}) is biholomorphic
here also. Moreover, since $\im A_k = z_{k+1}$ we see that (\ref{eq:Sigma}) maps $\tilde O_k$ to
$O_{k+1}$ and therefore $O_k$ maps to $\tilde O_{k+1}$. When $\lambda=0$ is a branch we choose
$\zeta=\sqrt\lambda$ to be a local parameter. A locally holomorphic  family $v_\zeta= v_0 +\zeta v_1 +
\ldots$ of eigenvectors now yields
\[
T_k^{\zeta^2} v_\zeta = -\zeta^{-2}A_kv_0 - \zeta^{-1}A_kv_1 + O(1).
\]
But $[v_0]=z_k$ since there is only one point over $\lambda=0$
hence $T_k$ has a simple pole at $O_k$. Again,
the image of $O_k$ under (\ref{eq:Sigma}) is $O_{k+1}$ since $z_{k+1}=\im A_k$.

Finally, since $T_k$ maps eigenlines to eigenlines it represents a rational section of
$\caL_k\otimes\caL_{k+1}^{-1}$ since $\caL_k$ is the dual of the eigenline bundle of $M_k$. By the
discussion above this section has divisor $S_k-\tilde O_k$ so 
$\caL_k\otimes\caL_{k+1}^{-1} \simeq\caO_\Sigma(S_k-\tilde O_k)$. But
(\ref{eq:Sigma}) maps $(O_k,\tilde O_k)$ to $(\tilde O_{k+1},O_{k+1})$ so we find that $O_k$ is
$O$ for $k$ even and $\tilde O$ for $k$ odd. This completes the proof.$\Box$

Let us define a periodic
map $L:\Z\to Jac(\Sigma)$ into the Jacobi variety (i.e.\ the group of isomorphism
classes of line bundles of degree zero over $\Sigma$) by setting $L_0=\caO_\Sigma$ and
$L_{k+1}\otimes L_k^{-1} = \caL_{k+1}\otimes\caL_k^{-1}$ (with $\caL_{k+n}=\caL_k$). The
previous proposition shows that when $n$ is odd this is a homomorphism, whereas when $n$ is even we
call it a {\em zigzag} since it is only a homomorphism on $2\Z$. In fact in either case
\begin{equation}
\label{eq:torsion}
\caO_\Sigma\simeq L_{2n}\simeq \caO_\Sigma(O-S+\tilde O-S)^n\simeq\caO_\Sigma(\tilde S-S)^n
\end{equation}
using the fact that $S+\tilde S\sim O+\tilde O$ (linear equivalence).
Therefore the divisor $\tilde S-S$ is a torsion divisor (in which case $S$ is called a division
point on $\Sigma$). Indeed $\tilde S-S$ satisfies a slightly stronger condition.
\begin{lem}
\label{lem:torsion}
The divisor $n(\tilde S -S)$ is the divisor of a rational function on $\Sigma$ which takes the same
value over the two points $P_\infty,\tilde P_\infty$ over $\infty$. 
\end{lem}
Another way of saying this is to say that $\tilde S-S$ is a torsion divisor on the 
singular curve $\Sigma^\prime$
obtained from $\Sigma$ by identifying the two points at infinity to obtain an ordinary double
point.

\noindent
{\bf Proof.} By the proof of the previous proposition $2S-(O+\tilde O)$ is the divisor of the rational
section of $\caL_{2j}\otimes\caL_{2j+2}^{-1}$ represented by $T_{2j+1}\circ T_{2j}$, 
therefore $S-\tilde S$ is the divisor of $(1-\lambda^{-1})^{-1}T_{2j+1}\circ T_{2j}$.
Observe that
\[
(H_0^\lambda)^2 = \prod_{j=n-1}^0 T_{2j+1}\circ T_{2j}
\]
(with the indices counted modulo $n$), so $(1-\lambda^{-1})^{-n}(H_0^\lambda)^2$ is a rational
section of $\caL_0\otimes\caL_0^{-1}$ with divisor $n(S-\tilde S)$. 
But $H_0^\lambda$ is itself a section of
$\caL_0\otimes\caL^{-1}_0\simeq\caO_\Sigma$ and from (\ref{eq:M}) we see that
any section $v$ of $\caL_0^{-1}$ satisfies
\[
H_0^\lambda v =(p+\lambda^{-1}\mu)v,
\]
where we recall that $p(\lambda)=\Tr(H_0^\lambda)/2$. Therefore
$H_0^\lambda$ represents $p+\lambda^{-1}\mu$, which takes the value $1$ at
any point where $\lambda^{-1}=0$ since $p(\infty)=1$. Therefore $n(S-\tilde
S)$ is the divisor for $(1-\lambda^{-1})^{-n}(p+\lambda^{-1}\mu)^2$, which also takes the value $1$
wherever $\lambda^{-1}=0$.$\Box$ 

\subsection{Recovery of the discrete curve from its spectral data.}

The spectral data $(\Sigma,\lambda,\caL,O,S)$ determines each 
$\caL_k$ by proposition \ref{pp:flow} if
we take $\caL_0=\caL$. This is enough to give each map $f_k:\Sigma\to\P^1$ up to a M\" obius transform:
we take it from the natural map $\Sigma\to\P\Gamma(\caL_k)^*$ which assigns to each point $P$ the
hyperplane $\Gamma(\caL_k(-P))$. 
By lemma \ref{lem:zk} we know $f_k(O_k)$ gives $z_k$ upon an appropriate identification of
$\P\Gamma(\caL_k)$ with $\P^1$. So to recover the curve $\Gamma$ we need only understand how this
identification is fixed. Indeed it is clear that since $\Gamma$ is only to be determined up to M\"
obius transformation what we really want to see is how each $\P\Gamma(\caL_k)$ is identified with,
say, $\P\Gamma(\caL_0)$. This is achieved by first 
identifying each space $\Gamma(\caL_k)$ with the sum of the two fibres
of $\caL_k$ over $\lambda=\infty$. We then interpret $T_k^\infty=I$ as identifying these 
fibres for $\caL_k$ with the fibres for $\caL_{k+1}$. 

More precisely, for each $k$ let $\caE_k\subset\Sigma\times\C^2$ denote the eigenline bundle with dual 
$\caL_k$. It follows that any linear form $e\in(\C^2)^*$, being a global section of 
$\Sigma\times(\C^2)^*$, induces a global section $\sigma_k$ of $\caL_k$.
Now let $\tau_k$ denote the section of ${\rm Hom}(\caE_k,\caE_{k+1})\simeq
\caL_k\otimes\caL_{k+1}^{-1}$ corresponding
to the map $T^\lambda_k$. Then $e\circ T^\lambda_k$ represents the rational section
$\sigma_{k+1}\otimes\tau_k$ of $\caL_k$. Since $T^\infty_k=I$ we have the
identity
\[
(\sigma_{k+1}\otimes\tau_k)\vert P = \sigma_{k}\vert P \quad \hbox{for $\lambda(P)=\infty$.}
\]
Here $\sigma\vert P$ denotes the section $\sigma$ restricted to $P$. This uniquely determines
$\sigma_{k+1}$ given $\sigma_{k}$ since no global section vanishes at both points over
$\lambda=\infty$ (by lemma \ref{lem:embedding}). Thus we have  maps
\begin{equation}
\label{eq:tau}
t_k:\Gamma(\caL_{k+1})\rightarrow\Gamma(\caL_k)\quad\hbox{where}\ t_k(\sigma)\vert\infty
=(\sigma\otimes\tau_k)\vert\infty
\end{equation}
and $\sigma\vert\infty$ denotes $(\sigma\vert P_\infty,\sigma\vert\tilde P_\infty)$.
This uses the identification 
\[
\Gamma(\caL_k)\rightarrow \caL_k\vert P_\infty\oplus\caL_k\vert \tilde P_\infty
\]
which restricts sections to the two fibres over infinity. 

For simplicity let $V_k$ denote the sum of
fibres on the right. Notice that to each point $P$ on 
$\Sigma$ we have a line in $V_k$, by evaluating the
section vanishing at $P$ at the two fibres over infinity. 
The lines for $P_\infty$ and $\tilde P_\infty$ are independent and we choose a third point $P$
to fix the identification of $\P V_0$ with $\P^1$ by sending these three lines to
$0,\infty$ and $1$ respectively. Combining this with the map $t_k$ from (\ref{eq:tau}) gives the 
identification of $\P
V_k$ with $\P^1$ for every $k$ (notice this only depends on the divisor of $\tau$ and not its
scale). 
Since $\caL$ has no sections which vanish at both $O,\tilde O$ any globally holomorphic
section of $\caL_k(-O_k)$ has divisor
$D_k+O_k$ where $D_k$ is a positive and non-special divisor of degree $g$.
\begin{lem}
\label{lem:Baker}
Given a discrete curve $\Gamma$ with spectral data as above,
let $\psi_k$ be the (unique up to scaling) non-zero
rational function on $\Sigma$ with divisor $D_k+ E_k-D_0$ where $E_k=\sum_{j=0}^{k-1}(S-O_j)$. 
Then we recover $\Gamma$, up to a M\"
obius transform, as the image of the map 
$z:\Z\to\P^1$ given by $z_k =\psi_k(P_\infty)/\psi_k(\tilde P_\infty)$.
\end{lem}
{\bf Proof.} Let $\sigma_k$ generate $\Gamma(\caL_k(-O_k))$, then $\sigma_k$ 
has divisor $D_k+O_k$. According
to (\ref{eq:tau}) it determines a line in $V_0$ by evaluating the section
$s_k=\sigma_k\otimes\tau_{k-1}\otimes\ldots\otimes\tau_0$ at $P_\infty$ and $\tilde P_\infty$. 
The resulting line $[s_k\vert
P_\infty,s_k\vert\tilde P_\infty]$ is then mapped to the line in $\P^1$ with homogeneous coordinates
$[(s_k/\sigma)\vert P_\infty, (s_k/\sigma)\vert\tilde P_\infty]$ where $\sigma$ is any section
generating the line $\Gamma(\caL_0(-P))$ corresponding to our third point $P$,
according to the identification of $\P V_0$ with $\P^1$ fixed above. Since $\psi_k = s_k/s_0
=(s_k/\sigma)/(s_0/\sigma)$ this rational function determines a 
M\" obius equivalent discrete curve. This function has divisor 
\[
D_k+O_k+\sum_{j=0}^{k-1}(S-\tilde O_j) - D_0-O_0 = D_k+\sum_{j=0}^{k-1}(S-O_j)-D_0,
\]
since $\tilde O_j=O_{j+1}$.$\Box$.

Notice that we will not have $\psi_k(P_\infty)=0=\psi_k(\tilde P_\infty)$ since we have assumed
that $\caL(-P_\infty-\tilde P_\infty)$ has  no global sections. 
We will postpone the explicit computation of the $z_k$ until we have 
introduced discrete conformal maps.

\section{Discrete Conformal Maps.}

A discrete conformal map is a map $z:\Z^2\to\P^1$ with the property that
\begin{equation}
\label{eq:dcm}
[z_{k,m+1}:z_{k,m}:z_{k+1,m}:z_{k+1,m+1}] = q
\end{equation}
for some constant $q\neq 0,1,\infty$ for all $k,m$.
We will be principally concerned with discrete conformal maps with one period i.e.\
we will assume there is an $n$ such that $z_{k+n,m} = z_{k,m}$ for all $k,m\in\Z^2$. In that case we
can also think of the map as describing the conformal flow of the discrete curve $\Gamma_{0,0} =
(z_{0,0},\ldots,z_{n-1,0})$. 

To each discrete curve $\Gamma_{k,m}$ in this flow let us assign its
spectral data $(\Sigma_{k,m},\lambda_{k,m},\caL_{k,m},O_{k,m},S_{k,m})$. 
\begin{lem} 
The point $Q_{k,m} = (q,z_{k,m+1})$ lies on $\Sigma_{k,m}$.
\end{lem}
{\bf Proof.} By (\ref{eq:dcm}) we see that
\[
T_{k,m}^q(z_{k,m+1}) = z_{k+1,m+1}
\]
for all $k,m$. It follows that $H_{k,m}^q(z_{k,m+1}) = z_{k,m+1}$.$\Box$

As earlier, we use $(\Sigma,\lambda)$ to denote the isomorphism class of $(\Sigma_{0,0},\lambda_{0,0})$.
We define $Q\in\Sigma$ to be the point corresponding to $Q_{0,0}$ on $\Sigma_{0,0}$. 
\begin{prop}
\label{pp:qflow}
For each $k,m$ there is an isomorphism $(\Sigma_{k,m},\lambda_{k,m})\simeq (\Sigma,\lambda)$ such
that $S_{k,m}$ is mapped to $S$, $Q_{k,m}$ is mapped to $Q$ 
but $O_{k,m}$ is mapped to $O$ for $k+m$ even and $\tilde O$ for
$k+m$ odd. Further:
\begin{equation}
\label{eq:L}
\begin{array}{c}
\caL_{k+1,m}\otimes\caL_{k,m}^{-1}\simeq \caO_\Sigma(\tilde O_{k,m}-S); \\
\caL_{k,m+1}\otimes\caL_{k,m}^{-1}\simeq \caO_\Sigma(\tilde O_{k,m}-Q). 
\end{array}
\end{equation}
\end{prop}
The proof of this is identical to the proof of proposition \ref{pp:flow}
given the next
lemma, which tells us how the holonomy changes under the conformal flow. Let us introduce the map
$\hat T_{k,m}:\P^1\to\P GL_2$ characterised by
\[
[z:z_{k,m}:z_{k,m+1}:\hat T^\lambda_{k,m}(z)]=\lambda.
\]
By earlier remarks this has matrix representation
\[
\hat T^\lambda_{k,m} =I-\lambda^{-1}\hat A_{k,m}
\]
where $\hat A_{k,m}$ is the projection matrix with kernel $z_{k,m}$ and image $z_{k,m+1}$. The following
lemma tells us how the holonomy evolves as we change the base point (cf. \cite{HerHP} for a similar result
about discrete isothermic nets).
\begin{lem}
The trace free part $M_{k,m}^\lambda$ of the holonomy for $\Gamma_{k,m}$ evolves according to
\begin{equation}
\label{eq:evolve}
\begin{array}{c}
M_{k+1,m}^\lambda = T_{k,m}^\lambda\circ M_{k,m}^\lambda\circ (T_{k,m}^\lambda)^{-1};\\
M_{k,m+1}^\lambda = \hat T_{k,m}^{\lambda/q}\circ M_{k,m}^\lambda\circ 
(\hat T_{k,m}^{\lambda/q})^{-1}. 
\end{array}
\end{equation}
\end{lem}
{\bf Proof.} The first identity we know from earlier. To prove the second identity it 
suffices to show that
\[
\hat T^{\lambda/q}_{k+1,m}\circ T^\lambda_{k,m} = T^\lambda_{k,m+1}\circ\hat T^{\lambda/q}_{k,m},
\]
when (\ref{eq:dcm}) holds. If we expand the matrix representations for these maps as functions of
$\lambda^{-1}$ we see that this is equivalent to showing that:
\[
\begin{array}{l}
\hbox{(a)}\ \hat A_{k+1,m}A_{k,m} = A_{k,m+1}\hat A_{k,m},\ \hbox{and}\\
\hbox{(b)}\ A_{k,m} + q\hat A_{k+1,m} = A_{k,m+1} + q\hat A_{k,m}.
\end{array}
\]
In (a) it is clear that on both sides the image of the first
matrix is the kernel of the second, hence both sides are identically zero. For (b) we
can compute the matrices explicitly. But this can be made easier by first mapping
$(z_{k,m+1},z_{k,m},z_{k+1,m},z_{k+1,m+1})$ to $(\infty,1,0,q)$ by M\" obius transform. If we lift
$z\in\P^1$ to $(z,1)^t$ (or $(1,0)^t$ when $z=\infty$), then elementary calculations show that:
\[
\begin{array}{cc}
A_{k,m} = \left(\begin{array}{cc}0&0\\ -1&1\end{array}\right);  &
q\hat A_{k+1,m} = \left(\begin{array}{cc}q&0\\ 1&0\end{array}\right); \\
A_{k,m+1} = \left(\begin{array}{cc}0&q\\ 0&1\end{array}\right) ; &
q\hat A_{k,m} = \left(\begin{array}{cc}q&-q\\ 0&0\end{array}\right) .
\end{array}
\]
The identity required follows immediately.
$\Box$

\smallskip\noindent
By arguments identical to those in the proof of proposition \ref{pp:flow} we see that $\hat
T^\lambda_{k,m}$ represents a rational section of $\caL_{k,m+1}^{-1}\otimes\caL_{k,m}$ with 
divisor $Q_{k,m}-\tilde O_{k,m}$ and combining this with proposition \ref{pp:flow} we deduce 
proposition \ref{pp:qflow}. Naturally this means the complete spectral data for a periodic discrete
conformal map is the sextuple $(\Sigma,\lambda,\caL,O,S,Q)$.
We will see later that any sextuple $(\Sigma,\lambda,\caL,O,S,Q)$ possessing
the properties of proposition \ref{pp:curve} is spectral data for a discrete conformal map.

\subsection{Explicit formula for the discrete conformal map.}

Given a discrete conformal map with generic spectral data we can write down an explicit formula for
it (up to M\" obius equivalence) in terms of the Riemann theta function of $\Sigma$. For this we
need the analogue of lemma \ref{lem:Baker} proved earlier. Recall we identify each
$\Gamma(\caL_{k,m})$ with $\Gamma(\caL_{0,0})$ in the following way. To $T_{k,m}^\lambda$ and 
$\hat T_{k,m}^\lambda$ we have
corresponding sections $\tau_{k,m}$ of $\caL_{k,m}\otimes\caL_{k+1,m}^{-1}$ and $\hat\tau_{k,m}$ of
$\caL_{k,m}\otimes\caL_{k,m+1}^{-1}$. Since every section $\sigma\in\Gamma(\caL_{k,m})$ is determined
entirely by its restriction $\sigma\vert\infty=(\sigma\vert P_\infty,\sigma\vert\tilde P_\infty)$ 
we may define bijective linear maps
\begin{equation}
\label{eq:tmap}
\begin{array}{cc}
t_{k,m}:\Gamma(\caL_{k+1,m}) \rightarrow  \Gamma(\caL_{k,m}) &  \hbox{where   
$t_{k,m}(\sigma)\vert\infty = (\sigma\otimes\tau_{k,m})\vert\infty$;}\\
\hat t_{k,m}:\Gamma(\caL_{k,m+1})  \rightarrow  \Gamma(\caL_{k,m}) & \hbox{where $\hat
t_{k,m}(\sigma)\vert\infty = (\sigma\otimes\hat\tau_{k,m})\vert\infty$.}
\end{array}
\end{equation}

To write this is terms of global
sections we need to introduce the following
divisors. First, let $D_{k,m}$ be the unique positive divisor in the linear system of
$\caL_{k,m}(-O_{k,m})$. Now define the
divisors $E_{k,m}$ by taking $E_{0,0}$ to be trivial and making
\[
\begin{array}{l}
E_{k+1,m}-E_{k,m} = S- O_{k,m}\\
E_{k,m+1}-E_{k,m} = Q- O_{k,m}
\end{array}
\]
\begin{lem}
\label{lem:qBaker}
Up to M\" obius equivalence the 
discrete conformal map with spectral data
$(\Sigma,\lambda,\caL,O,S,Q)$ is given by the map $z:\Z^2\to\P^1$ for which
$z_{k,m} = \psi_{k,m}(P_\infty)/\psi_{k,m}(\tilde P_\infty)$, where $\psi_{k,m}$
is the (unique up to scaling) rational function on $\Sigma$ with divisor 
\begin{equation}
\label{eq:divisor}
(\psi_{k,m}) = D_{k,m} +E_{k,m}  - D_{0,0}.
\end{equation}
\end{lem}
The proof is the same as for lemma \ref{lem:Baker}.

To explicitly compute the $z_{k,m}$, we need to
fix a basis $\{a_j,b_j\}_{j=1}^g$ for the first homology of
$\Sigma$ with the standard intersection pairing. We use the a-cycles to fix a
dual basis $\{\omega_j\}_{j=1}^g$ of holomorphic one forms and thereby equate
$Jac(\Sigma)$ with $\C^g/\Lambda$ where $\Lambda$ is the lattice representing
the first homology via integration of the vector $(\omega_1,\ldots,\omega_g)$ over each
homology class. Given a base point $B$ on $\Sigma$ we have the Abel map
\[
\caA:\Sigma\rightarrow \C^g/\Lambda\ ,\quad P\mapsto \int_B^P(\omega_1,\ldots,\omega_g)\bmod\Lambda,
\]
and, more generally, its extension to divisors by addition i.e.\
$\caA(P+Q)=\caA(P)+\caA(Q)$. Further, let $\omega_{g+1}$ be the
unique meromorphic differential satisfying:
(1) $\omega_{g+1}$ is holomorphic except at $P_\infty$ and $\tilde P_\infty$ 
where it has simple poles of
residue $1/2\pi i$ and $-1/2\pi i$ respectively;
(2) its integral around any a-cycle is zero. 

Given all this, we define maps $\alpha_j:\Z^2\to\C$ for $j=1,\ldots,g+1$ by setting $\alpha_j(0,0) =0$
and 
\[
\begin{array}{c} 
\alpha_j(k+1,m)-\alpha_j(k,m)= \int_{O_{k,m}}^S\omega_j;\\
\alpha_j(k,m+1)-\alpha_j(k,m)= \int_{O_{k,m}}^Q\omega_j.
\end{array}
\]
In the formula to follow we will write $\alpha^\prime_{k,m}:\Z^2\to\C^{g+1}$ for the map whose $j$-th
component is $\alpha_j(k,m)$, for $j=1,\ldots,g+1$ and $\alpha_{k,m}:\Z^2\to\C^g$ 
for its projection onto the first $g$ components. We think of the former as lying over the generalised
Jacobian $J^\prime$ of the curve $\Sigma^\prime$ obtained by identifying 
$P_\infty$ with $\tilde P_\infty$ to obtain a node. The group $J^\prime$ may be analytically realised
as $\C^{g+1}/\Lambda^\prime$ where $\Lambda^\prime$ represents the first homology of
$\Sigma-\{P_\infty,\tilde P_\infty\}$, the open variety of smooth points on $\Sigma^\prime$, 
via integration of the augmented vector $(\omega_1,\ldots,\omega_{g+1})$ (see e.g.\ \cite{Ser} p101). 
This point of view is
useful for computing the periodicity conditions.
\begin{thm}
\label{th:formula}
The formula
\begin{equation}
\label{eq:theta}
z_{k,m} = \exp[2\pi i\alpha_{g+1}(k,m)]\frac{\theta(\caA(P_\infty)+\alpha_{k,m}-
\caA(D_{0,0}) - \kappa)}{\theta(\caA(\tilde P_\infty)+\alpha_{k,m}-\caA(D_{0,0})-\kappa)},
\end{equation}
where $\kappa$ is the vector of Riemann constants,
recovers the discrete conformal map with spectral data $(\Sigma,\lambda,\caL,O,S,Q)$ up to M\" obius
transform. This map is periodic (in $k$) with period $n$ precisely when
$\alpha^\prime_{k+n,m}\equiv\alpha^\prime_{k,m}\bmod\Lambda^\prime$ for some (and hence all) $k,m$.
\end{thm}
{\bf Proof.} We begin by computing a function $\psi_{k,m}$ with divisor (\ref{eq:divisor}). For any
$A\in\Sigma$ let
$\eta^A_{k,m}$ be the unique meromorphic differential on $\Sigma$ with zero
a-periods and simple poles only at $A$ and $O_{k,m}$, where it has residues $1/2\pi i$
and $-1/2\pi i$ respectively. On the universal cover of $\Sigma$ we may define functions
$\beta_{k,m}$ by setting $\beta_{0,0}\equiv 0$ and 
\[
\begin{array}{ccc}
\beta_{k+1,m}(P)-\beta_{k,m}(P) & = & \int_B^P\eta^S_{k,m},\\ 
\beta_{k,m+1}(P)-\beta_{k,m}(P) & = & \int_B^P\eta^Q_{k,m}.
\end{array}
\]
Then define
\begin{equation}
\label{eq:Baker}
\psi_{k,m} = \exp(2\pi i\beta_{k,m}(P))
\frac{\theta(\caA(P)+\alpha_{k,m}-\caA(D_{0,0})-\kappa)}
{\theta(\caA(P)-\caA(D_{0,0})-\kappa)}
\end{equation}
where every integral from $B$ to $P$ is along the same path. By a standard reciprocity formula
for differentials (see e.g.\ \cite{GriH})
\[
\oint_{b_j}\eta^A_{k,m} = \int_{O_{k,m}}^A\omega_j .
\]
Using this it is a simple exercise to check that $\psi_{k,m}$ is well-defined on $\Sigma$. By
construction, the denominator vanishes precisely on $D_{0,0}$ and the 
exponential term contributes the divisor $E_{k,m}$. Moreover, by definition
$\alpha_{k,m} = \caA(E_{k,m})$ and $\caA(D_{k,m})=\caA(D_{0,0})-\caA(E_{k,m})$ so the numerator
vanishes precisely on $D_{k,m}$. Hence $\psi_{k,m}$ has divisor (\ref{eq:divisor}). By lemma
\ref{lem:qBaker} the discrete conformal map is, up to M\" obius transform, given by
$\psi_{k,m}(P_\infty)/\psi_{k,m}(\tilde P_\infty)$. However, by another reciprocity formula
\[
\int_{O_{k,m}}^A \omega_{g+1} = \int_{\tilde P_\infty}^{P_\infty}\eta^A_{k,m}
\]
hence $\beta_{k,m}(P_\infty)-\beta_{k,m}(\tilde P_\infty) = \alpha_{g+1}(k,m)$. Hence, up to
multiplication by a constant independent of $k,m$, we obtain (\ref{eq:theta}). 

For the periodicity, we will give a proof which works even when $\Sigma$ is singular, since we will
compute examples with such spectral curves shortly. First recall (from e.g.\ \cite{Ser})
that the curve $\Sigma^\prime$ has its own Abel map:
\[
\caA^\prime:\Sigma-\{P_\infty,\tilde P_\infty\}\rightarrow J^\prime,
\qquad P\mapsto \int_B^P(\omega_1,\ldots,\omega_{g+1})\bmod\Lambda^\prime.
\]
With this notation we have
\[
\alpha^\prime_{k+n,m}-\alpha^\prime_{k,m}  =  \caA^\prime(E_{k+n,m}-E_{k,m})
                             =\left\{
\begin{array}{l} \caA^\prime(nS-\frac{n}{2}O-\frac{n}{2}\tilde O), \hbox{ $n$ even;}\\
                 \caA^\prime(nS-nO), \hbox{ $n$ odd.}\end{array}\right.
\]
In particular this is independent of $k,m$. Also recall that Abel's theorem holds for $\Sigma^\prime$ i.e.\
$\caA^\prime(D)\equiv 0$ if and only if $D$ is the divisor of a rational function on $\Sigma^\prime$
(equally, $D$ is the divisor of a rational function on $\Sigma$ taking the same value at $P_\infty,\tilde
P_\infty$). Now the rational function 
$\psi_{k+n,m}/\psi_{k,m}$ has divisor $D_{k+n,m}-D_{k,m}+E_{k+n,m}-E_{k,m}$ so if
$\caA^\prime(E_{k+n,m}-E_{k,m})\equiv 0$ then $D_{k+n,m}-D_{k,m}$ is itself the divisor of a rational
function on $\Sigma$. But each $D_{k,m}$ is positive of degree $g$ and non-special, therefore
$D_{k+n,m}=D_{k,m}$. Thus $\psi_{k+n,m}/\psi_{k,m}$ must take the same value at $P_\infty,\tilde
P_\infty$, whence $z_{k+n,m}=z_{k,m}$. On the other hand, we have seen earlier that when the map is
periodic the function $p+\lambda^{-1}\mu$ (representing the eigenvalues of the holonomy matrix) has
divisor $E_{k+n,m}-E_{k,m}$. Since $H^\infty_0=I$ this function takes the same values at $P_\infty,\tilde
P_\infty$.
$\Box$

\smallskip\noindent
{\bf Singularities and the $\theta$-divisor.} The formula (\ref{eq:theta}) need not give a discrete 
conformal map for all $k,m$. Indeed we expect there to be singularities when both translates of the
$\theta$-function are zero: this will occur on some codimension two subvariety of $\C^g$. However, and
perhaps less obviously, the map will also have singularities whenever
$\caL_{k,m}(-P_\infty-\tilde P_\infty)$ fails to be non-special i.e.\ on some translate
of the $\theta$-divisor. It is interesting to see what happens in this circumstance: typically
the map fails to be
a discrete immersion i.e.\ adjacent points fail to be distinct (see Figure 3). 

\begin{figure}[ht]
\centering
\includegraphics[scale=0.6]{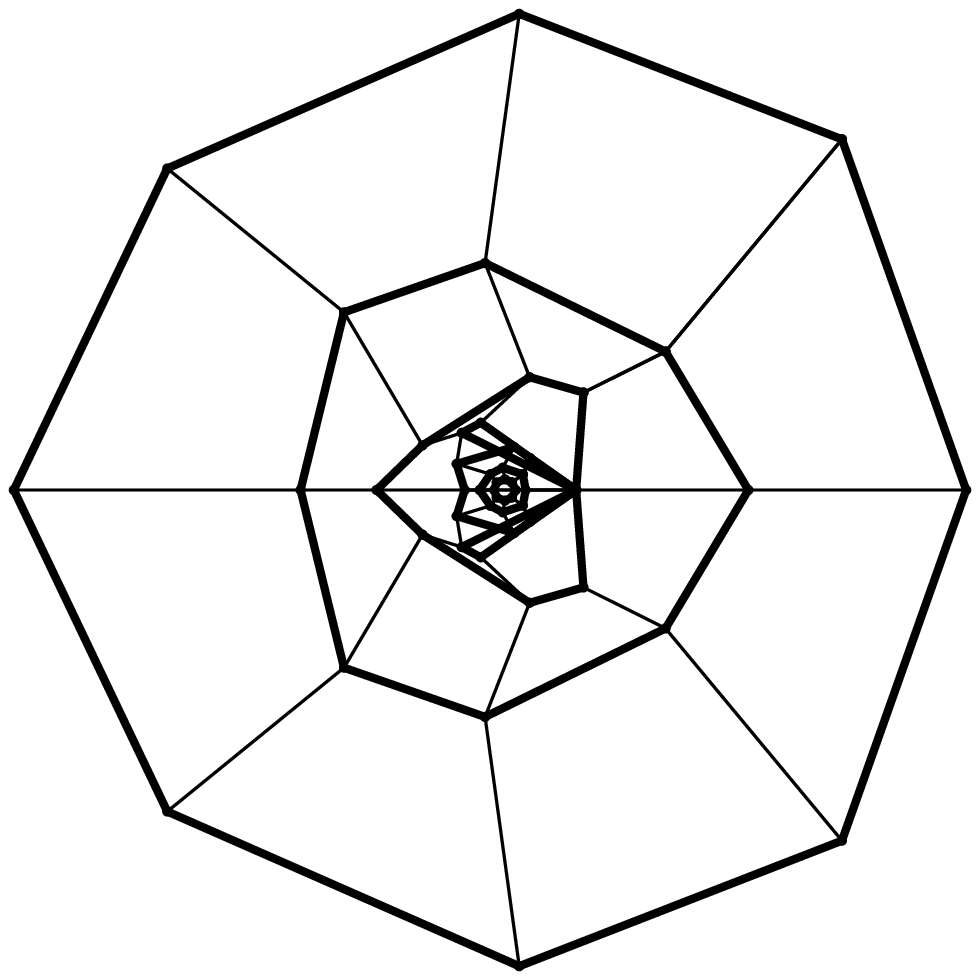}
\includegraphics[scale=0.5]{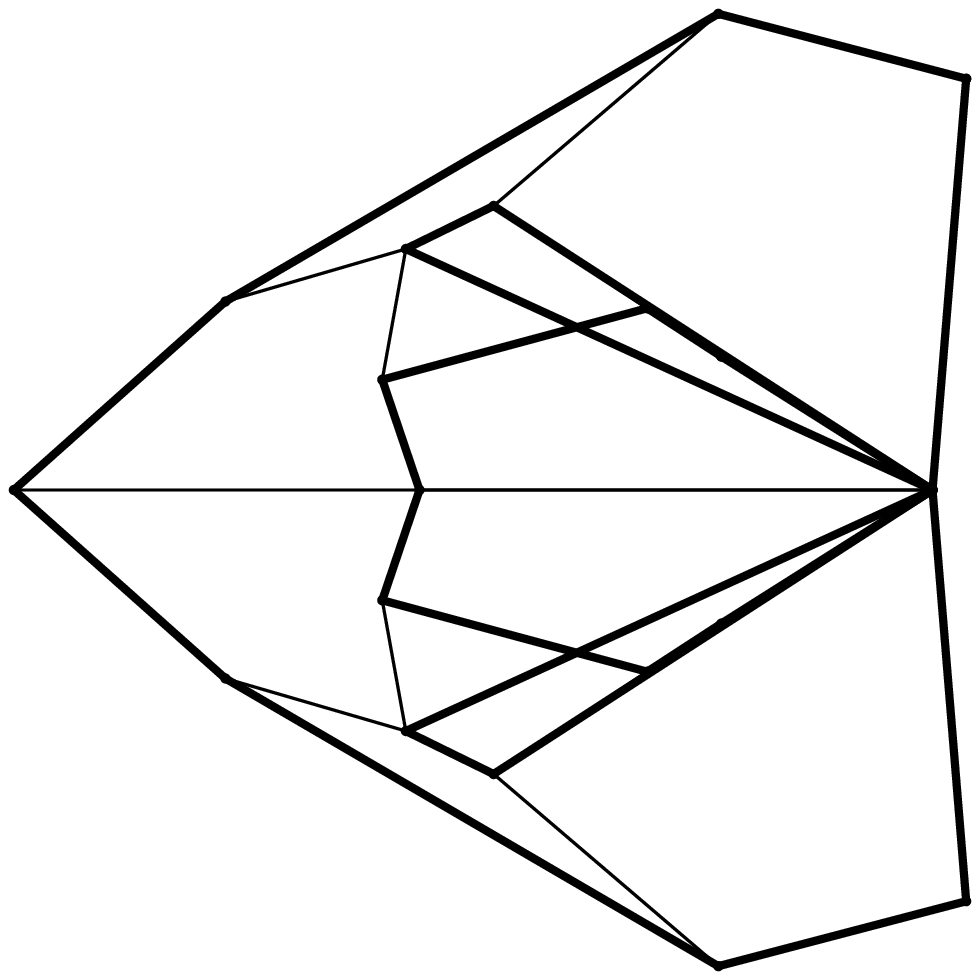}
\caption{The collapsing of points at a singularity (with close up on the right). Points with the same $m$ are joined
by bolder lines.}
\end{figure}

\noindent
Let us see why this occurs.
Suppose that $\caL(-P_\infty-\tilde P_\infty)$ is special  Since all pairs
$P+\tilde P$ are linearly equivalent this means  the divisor class for $\caL(-O-\tilde O)$ contains a
positive divisor (of degree $g-1$), $E$ say, and we take $D=E+\tilde O$. Now let us define
\[
f(P) = \exp(2\pi i\int_O^P\omega_{g+1})\frac{\theta(\caA(P)-\caA(D-P_\infty)-\kappa)}
                                {\theta(\caA(P)-\caA(D-\tilde P_\infty)-\kappa)}
\]
with the base point for the Abel map at $O$. In that case, unless $D-P_\infty,D-\tilde P_\infty$ are
both special (which is the codimension two condition mentioned above),
this is a well-defined rational function on
$\Sigma$ and a careful comparison with (\ref{eq:theta}) shows that 
\[
z_{0,0}=f(O),\ z_{1,0}=f(S),\ z_{-1,0}=f(\tilde S),\ z_{0,1}=f(Q),\ z_{0,-1} = f(\tilde Q).
\]
The divisor for $f$ is $P_\infty -\tilde P_\infty +C-C^\prime$ where $C,C^\prime$ are the unique positive
divisors of degree $g$ for which $C-O\sim D-P_\infty$ and $C^\prime-O\sim D-\tilde P_\infty$. But
$D=E+\tilde O$, $O+\tilde O\sim P_\infty +\tilde P_\infty$ and $C,C^\prime$ are unique, hence $C=E+\tilde
P_\infty$ and $C^\prime =E+P_\infty$. Therefore the poles and zeroes of $f$ cancel i.e.\ $f$ is constant, 
so the five points above are identical.
In that case the cross-ratio condition breaks down locally.

\medskip\noindent
{\bf Geometric interpretation.} One of the reasons for choosing to write (\ref{eq:theta}) in 
this form is to exhibit an elegant geometric interpretation which supports the claim (made in 
\cite{NijC}) that the discrete conformal map equations are a discretization of the Schwarzian KdV (SKdV)
equations:
\begin{equation}
\label{eq:SKdV}
z_t=S(z)z_x, \qquad S(z)=z_{xxx}z_x^{-1}-\frac{3}{2}(z_{xx}z_x^{-1})^2.
\end{equation}
Here $S(z)$ is the Schwarzian derivative: it is well-known that $u(x,t)=2S(g)$ satisfies the KdV equation. 
The `finite gap' solutions of (\ref{eq:SKdV}) are related to
the formula (\ref{eq:theta}) in the following way. 

Let us assume that in the spectral data the point
$O$ is ramified (i.e.\ $O=\tilde O$). 
It is not hard to see that the projection $p:\C^{g+1}\to\C^g$ onto
the first $g$ coordinates projects the lattice $\Lambda^\prime$ onto $\Lambda$ and therefore it
induces a surjective homomorphism $p:J^\prime\to Jac(\Sigma)$ whose kernel
is $\Lambda^\prime/\Lambda\simeq\C/\Z\simeq \C^*$. 
Further, one can compute by means of multipliers that the 
pullback of the $\theta$-line bundle over $Jac(\Sigma)$ to $J^\prime$ has its space of globally
holomorphic sections spanned by the infinite collection
\[
\{\theta_k(Z)=\exp(2\pi ikZ_{g+1})\theta(p(Z) +k\caA(P_\infty-\tilde P_\infty)) :\ k\in\Z,\
Z\in\C^{g+1}\}.
\]
We see, therefore, that if we let $U_P$ denote $\caA^\prime(P-O)$ then up to a scaling the 
formula (\ref{eq:theta}) can be re-expressed as
\[
z_{k,m} =\frac{\theta_1(kU_S + mU_Q +\tau)}
         {\theta_0(kU_S + mU_Q +\tau)}
\]
for some constant $\tau\in\C^{g+1}$. Geometrically $U_S$ and $U_Q$ are the secants $\vec{OS}$ and
$\vec{OQ}$ on $\caA^\prime(\Sigma_0)$ respectively. 

On the other hand, if $U_1,U_3\in\C^{g+1}$ represent, respectively, the tangent 
to $\caA^\prime(\Sigma_0)$ at $O$ and its third derivative there (i.e.\
$U_1=(\partial\caA^\prime/\partial\zeta)(0)$ and $U_3=(\partial^3\caA^\prime/\partial\zeta^3)(0)$ for the 
local parameter $\zeta=\sqrt\lambda$), then
\[
z(x,t) =\frac{\theta_1(xU_1 + tU_3 +\tau)}
         {\theta_0(xU_1 + tU_3 +\tau)}
\]
satisfies the SKdV equation (\ref{eq:SKdV}) (cf.\ the formula given at the end of
\cite{McI}). Hence when viewed in the Jacobi variety (that is to say, after the equations have been
linearised) the discretization is nothing other 
than the perturbation of a tangent into a secant on the spectral curve. 

\subsection{Examples.}

Here we will present three examples, all based on taking $\Sigma$ to be a rational nodal
curve. The theory
works equally well in this case and it is a good deal easier to calculate since
the $\theta$-functions are simply polynomials in exponentials. Also, the 
periodicity condition is easy to satisfy and we can obtain discrete maps with any period. 

\noindent{\bf 1. $\Sigma$ is the Riemann sphere.} Take $\Sigma$ to be the Riemann sphere equipped with
an affine coordinate $\zeta$. In this case the $\theta$-function is a constant which we may as well
take to equal $1$. We take the hyperelliptic involution to be $\zeta\mapsto -\zeta$ and prescribe the
spectral data by $\zeta(O)=\epsilon$, $\zeta(S)=a$, $\zeta(Q)=b$ and $\zeta(P_\infty)=y$ where these
are all distinct and different from $0,\infty$. It follows
that
\[
\lambda = \frac{(a^2-y^2)}{(\zeta^2-y^2)}\frac{(\zeta^2-\epsilon^2)}{(a^2-\epsilon^2)}.
\]
Since $g=0$ in this case we only need to compute $\omega_1$ and its integrals. It is easy to see
that $\omega_1=\omega_y$ where
\begin{equation}
\label{eq:oneform}
\omega_y = \frac{1}{2\pi i}(\frac{1}{\zeta-y}-\frac{1}{\zeta + y})d\zeta
\end{equation}
and therefore, following the procedure above, we see that $z_{k,m}=h_{k,m}(y)$ where
\begin{equation}
\label{eq:hom}
h_{k,m}(y) = \exp[2\pi i\alpha_1(k,m)] =\left\{\begin{array}{cc}
(\frac{a-y}{a+y})^k(\frac{b-y}{b+y})^m(\frac{\epsilon +y}{\epsilon -y}) & k+m\ \hbox{odd},\\
(\frac{a-y}{a+y})^k(\frac{b-y}{b+y})^m & k+m\ \hbox{even}.
\end{array}
\right.
\end{equation}
This map has cross-ratio $\lambda(Q)$ 
and period $n$ whenever $(a+y)/(a-y)$ is an $n$-th root of unity ($n$ must be even if
$\epsilon\neq\infty$). Observe that this agrees
with the periodicity condition given in theorem 
\ref{th:formula}, for in this example $\Sigma^\prime$ is a
one node curve and  $J^\prime=\C/\Z\langle\oint\omega_1\rangle$, where the integral is around any
cycle separating $y$ from $-y$. It is not hard to see that parameter values can be chosen 
to give any cross-ratio with any period
$\geq 4$. 

This is the general case of the zigzag: when $\epsilon =\infty$ we have a homomorphism $\Z^2\to\C^*$. 
It is clear that this
is the discrete exponential function $k,m\mapsto \exp(kA + mB)$ for constants
$A,B$. 

\begin{figure}[ht]
\centering
\includegraphics[scale=0.5]{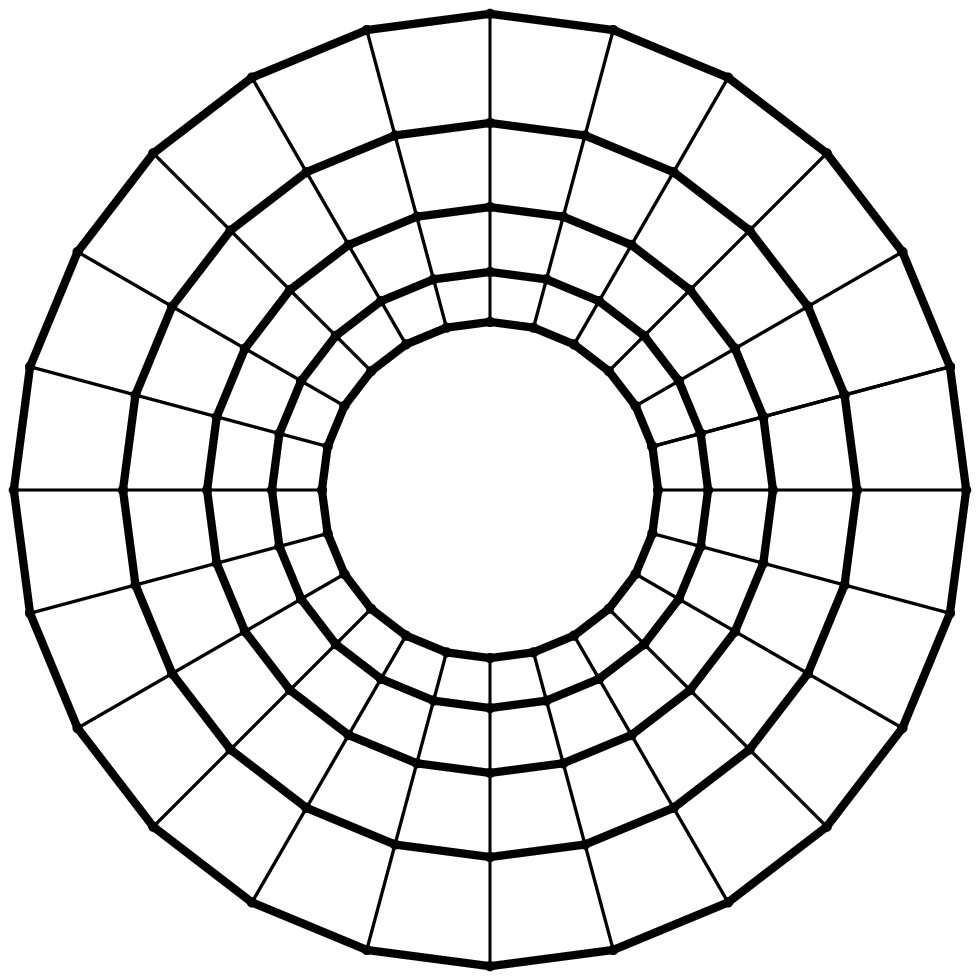}
\includegraphics[scale=0.5]{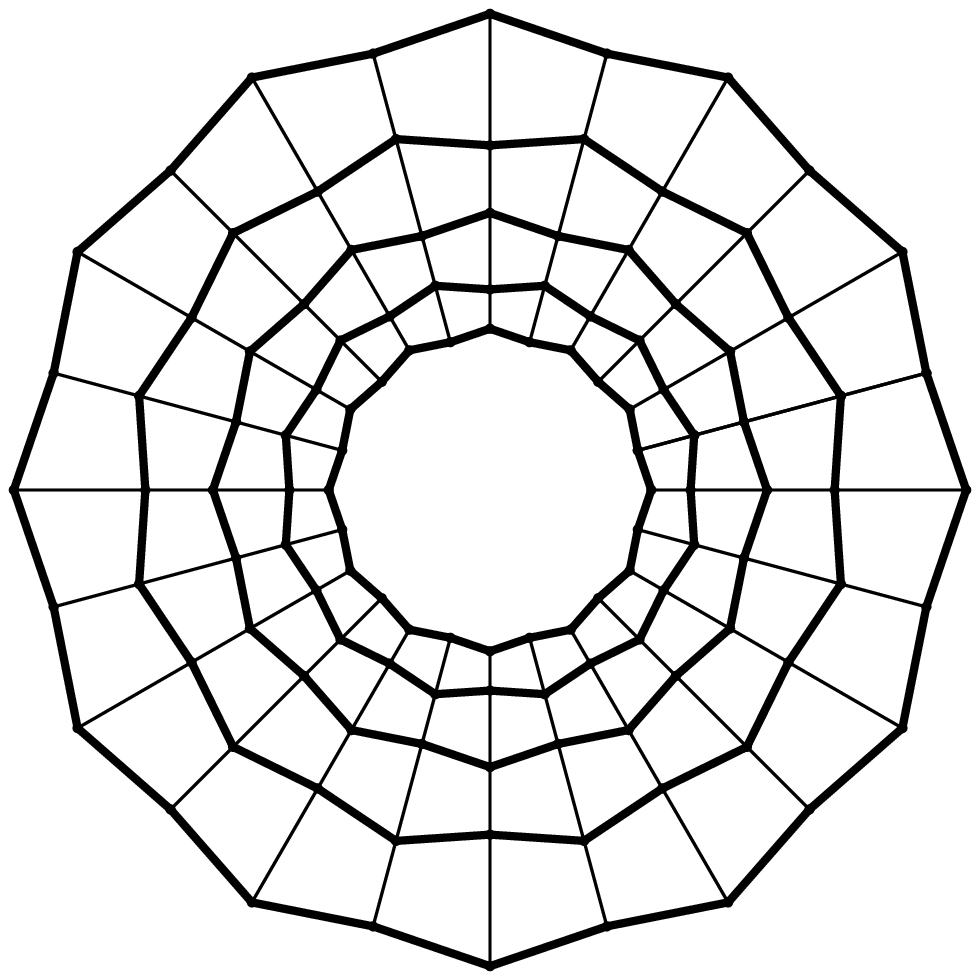}
\caption{The discrete exponential (left) and a zigzag version (right).}
\end{figure}

\medskip\noindent
{\bf 2. $\Sigma$ is a rational curve with one node.} Here we take $\Sigma$ to be $\P^1$ with the points
$\pm x\neq 0,\infty$ identified together to obtain a node. In this case
the $\theta$-function is given by $\theta(Z)=\exp(2\pi iZ)-1$. Here the space of regular one forms
corresponds to the meromorphic forms on $\P^1$ with simple poles at $x,-x$ only (see \cite{Ser} p68), 
hence the space is one dimensional and the (arithmetic) genus is $g=1$. 
We keep the same notation as in the previous
example and make the same choice for $\lambda$ (so $a,b,x,y,\epsilon$ are all distinct). 
One easily checks that the appropriate choices
for $\omega_1,\omega_2$ in this case are to take $\omega_1=\omega_x$ and $\omega_2=\omega_y$ using
(\ref{eq:oneform}). By computing the integrals and applying the formula (\ref{eq:theta}) we see that we
can write 
\[
z_{k,m} = h_{k,m}(y)\left(\frac{\left(\frac{y-x}{y+x}\right)h_{k,m}(x)-e^{2\pi ic}}
{\left(\frac{y+x}{y-x}\right)h_{k,m}(x)- e^{2\pi ic}}\right)
\]
where $h_{k,m}(y)$ is given by (\ref{eq:hom}) and the constant $c$ represents the terms
$\caA(D_{0,0})+\kappa$ in (\ref{eq:theta}). In this case $\Sigma^\prime$ is a two node curve 
with Jacobian 
$J^\prime=\C/\Z\langle \oint\omega_1,\oint\omega_2\rangle$ where the integrals are around $x$ and $y$
respectively. According to theorem \ref{th:formula} the map is periodic when both $(a-y)/(a+y)$ and
$(a-x)/(a+x)$ are $n$-th roots of unity (distinct, so that $x\neq\pm y$). In this case the
periodicity problem can be solved for any period $\geq 5$ for any value of the cross-ratio $\lambda(Q)$. 

\begin{figure}[ht]
\centering
\includegraphics[scale=0.3]{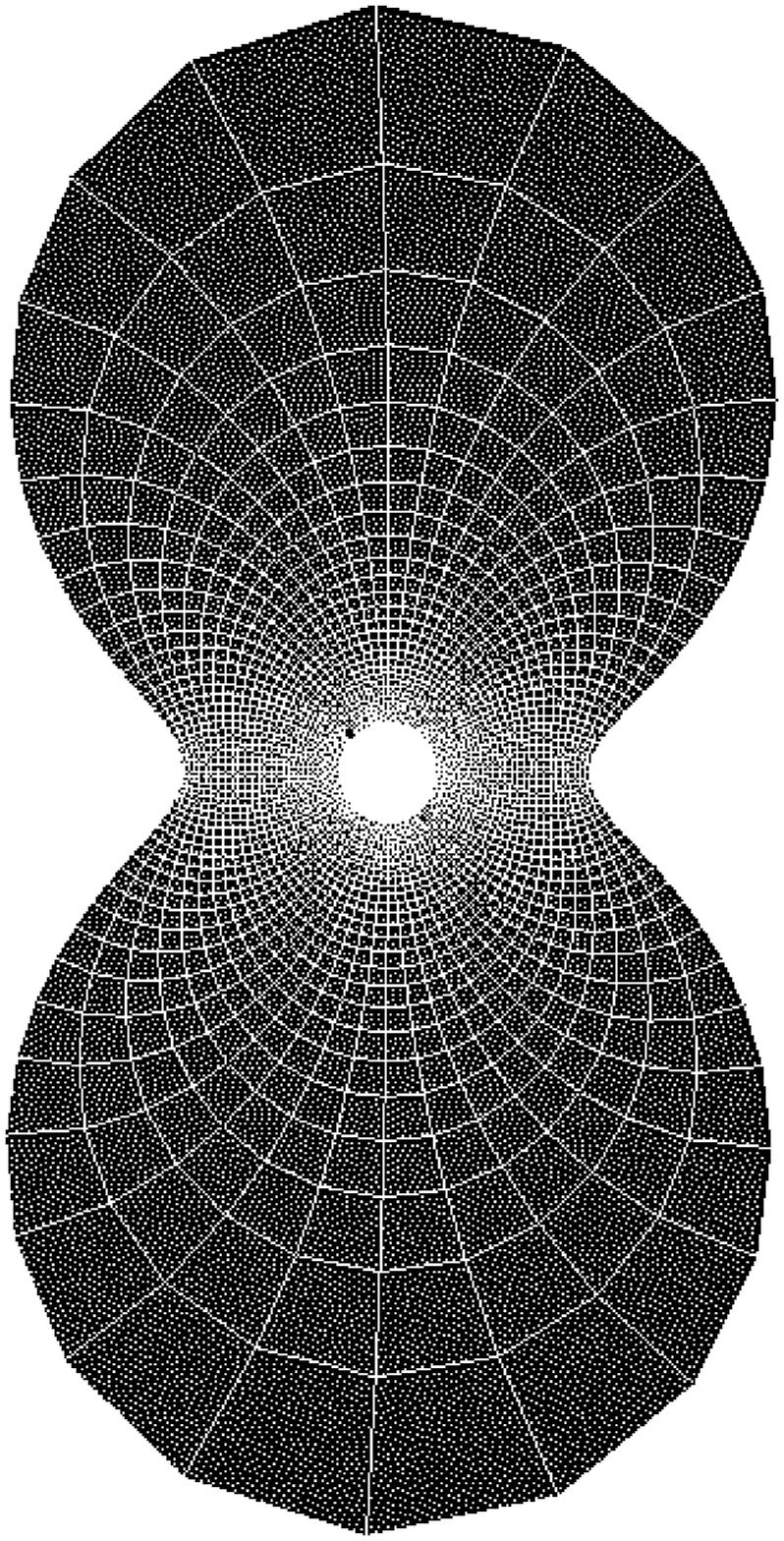}
\includegraphics[scale=0.3]{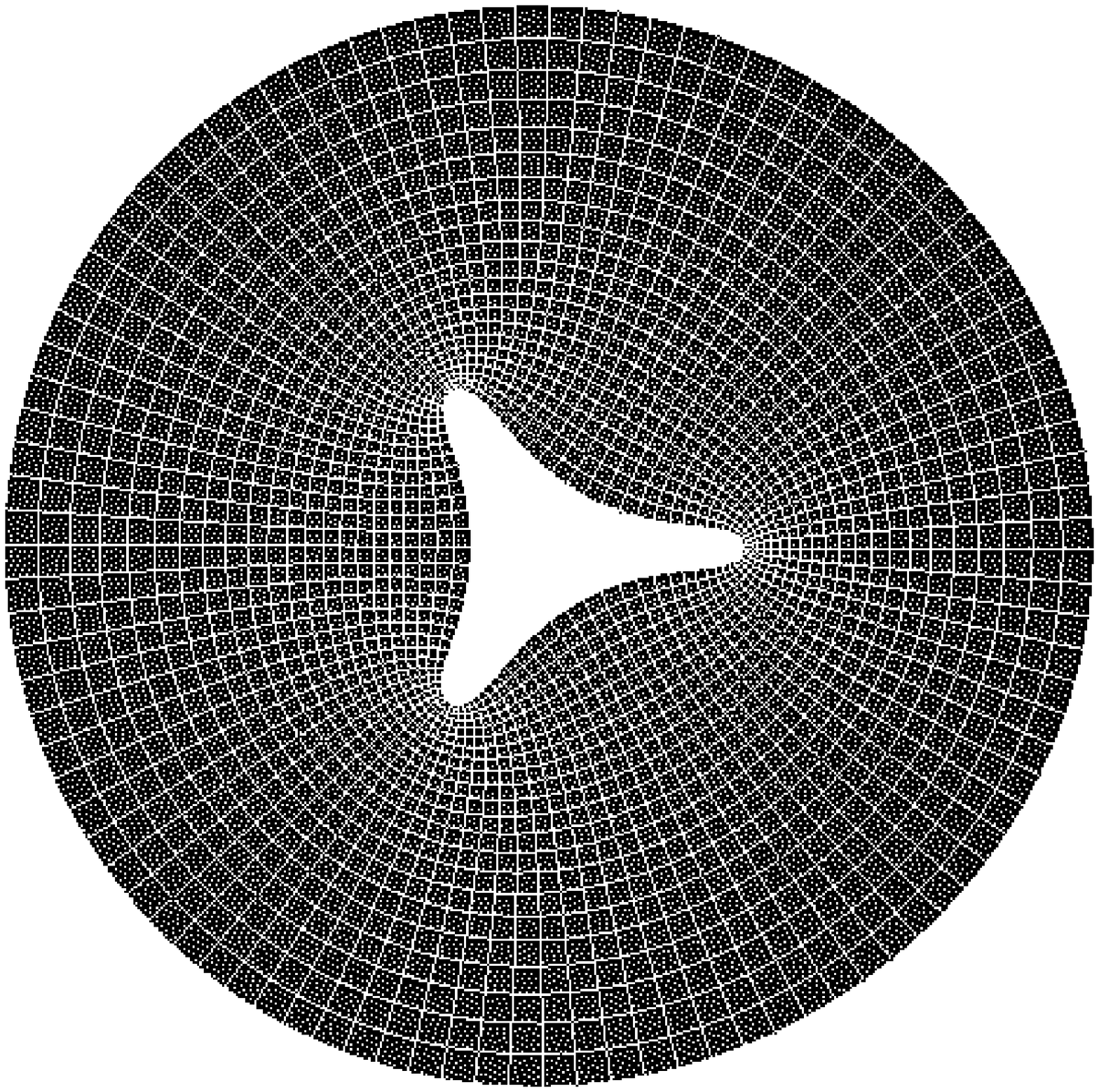}
\caption{Two discrete 1-solitons: showing 2-fold and 3-fold symmetries.}
\end{figure}

Each of these maps behaves like a discrete 1-soliton in the sense that it has 
asymptotics in $m$ like the discrete
exponential (see Figure 5). Let $B =(b-y)/(b+y)$, then for $\vert\ B\vert <1$ we have
\[
z_{k,m}/h_{k,m}(y)\rightarrow\left\{\begin{array}{l}
1\ \hbox{as $m\to\infty$;}\\
(y-x)^2/(y+x)^2\ \hbox{as $m\to-\infty$.} \end{array}\right.
\]
For $\vert B\vert >1$ the limits are interchanged. 

\begin{figure}[ht]
\centering
\includegraphics[scale=0.4]{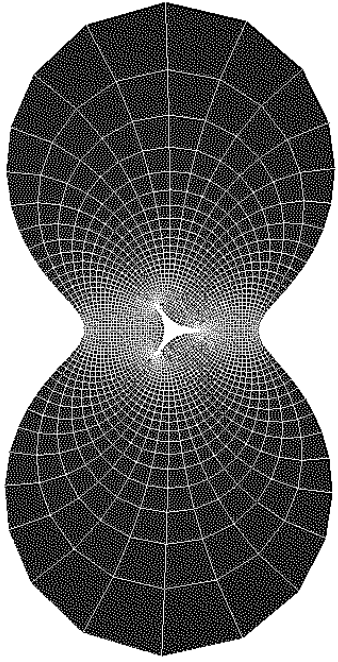}
\caption{A discrete 2-soliton: superposition of 2-fold and 3-fold symmetries.}
\end{figure}

\medskip\noindent
{\bf 3. $\Sigma$ is a rational curve with two nodes.} Take $\Sigma$ to be $\P^1$ with $\pm x_1$, $\pm x_2$
identified in pairs. The $\theta$-function here is given by 
\[
\theta(Z)=F(e^{2\pi iZ_1},e^{2\pi iZ_2}),\qquad F(X,Y)=\det\left(\begin{array}{cc} X-1 & (Y+1)x_1 \\
Y-1 & (X+1)x_2\end{array}\right).
\]
The arithmetic genus is $g=2$ and $(\omega_1,\omega_2,\omega_3)=(\omega_{x_1},\omega_{x_2},\omega_y)$ using
(\ref{eq:oneform}). Again, with $\lambda$ chosen as above the appropriate computation yields
\[
z_{k,m} = h_{k,m}(y)\frac{F(
(\frac{y-x_1}{y+x_1})h_{k,m}(x_1)e^{-2\pi ic_1},
(\frac{y-x_2}{y+x_2})h_{k,m}(x_2)e^{-2\pi ic_2})}
{F( (\frac{y+x_1}{y-x_1})h_{k,m}(x_1)e^{-2\pi ic_1},
(\frac{y+x_2}{y-x_2})h_{k,m}(x_2)e^{-2\pi ic_2})}
\]
where $c_1,c_2$ are parameters corresponding to the initial point $\caA(D_{0,0})+\kappa$. The
periodicity conditions can be solved for any period $\geq 7$ and any cross-ratio $\lambda(Q)$. 
Each of these maps behaves like a 2-soliton in the sense that: (a)
as $m\to\pm\infty$ it behaves like the discrete exponential, and; (b) by suitable choice of $c_1,c_2$
we observe that it behaves like two interacting 1-solitons (see Figure 6).

\subsection{Spectral data which produces discrete conformal maps.}

The aim of this section is to prove that any choice of spectral data 
satisfying the conditions described above yields a periodic discrete conformal map. 
\begin{thm}
Let $\Sigma$ be a compact hyperelliptic Riemann surface of genus $g$ with: 
a degree two function $\lambda$ unbranched at
$\lambda=1,\infty$; a degree $g+1$ line bundle $\caL$ for which
$\caL(-P_\infty-\tilde P_\infty)$ is non-special;
and non-singular points $O,S,Q$ with $\lambda(O)=0$, $\lambda(S)=1$ and $\lambda(Q)\neq 0,1,\infty$.
Then the data $(\Sigma,\lambda,\caL,O,S,Q)$ is the spectral data for a discrete conformal map
$z:\Z^2\to\P^1$ of cross-ratio $\lambda(Q)$. This map is periodic with period $n$ if and only 
if either: a)  n is even and 
$\frac{n}{2}(2S-O-\tilde O)$ is the divisor of a 
rational function on $\Sigma$ taking the same value at each point
over $\lambda=\infty$; or, b) n is odd and $n(S-O)$
is such a divisor. 
\end{thm}
{\bf Remarks.}(i) Strictly speaking it may be that $z_{k,m}$ has some singularities. These will occur
only  when $\caL_{k,m}(-P_\infty-\tilde P_\infty)$ (where $\caL_{k,m}$ is given by (\ref{eq:L})) is special i.e.\
has a non-zero global section. 

\noindent 
(ii) As we have seen in the examples above, the condition that $\Sigma$ be smooth can be weakened to 
allow irreducible singular algebraic curves, provided the points $O,\tilde O,S,Q$ are all smooth. 
The proof we give works in this generality, in particular, we will not rely on the $\theta$-function
formula given earlier.
 
\smallskip
We can immediately deduce from this theorem a simple fact about smooth families of periodic discrete
conformal maps.
\begin{cor} 
\label{cor:family}
Set $\Sigma^o=\Sigma-\lambda^{-1}(\{0,1,\infty\})$. 
Up to M\" obius equivalence, each discrete conformal map lies in a g+1-parameter 
family $z(Q,\caL)$ parameterised by $\Sigma^o\times Jac(\Sigma)$. 
\end{cor}
Notice that if we swap $Q$ with its involute
$\tilde Q$ we get $z_{k,m}(\tilde Q,\caL) = z_{k,-m}(Q,\caL)$, since 
$\caO_\Sigma(\tilde Q-O)\simeq\caO_\Sigma(O-Q)$. 

\smallskip
Now let us turn to proving the theorem.
We saw in the previous sections that the point $z_{k,m}$ is the image of the hyperplane
$\Gamma(\caL_{k,m}(-O_{k,m}))\subset\Gamma(\caL_{k,m})$, where $O_{k,m}$ is $O$ for 
$k+m$ even and $\tilde O$ otherwise,
under an identification 
\begin{equation}
\label{eq:zeta}
\zeta_{k,m}:\P\Gamma(\caL_{k,m})^*\rightarrow\P^1
\end{equation}
of $\P^1$ with the space $\P\Gamma(\caL_{k,m})^*$ of hyperplanes 
in the two dimensional space $\Gamma(\caL_{k,m})$.
This identification is defined as follows. For each $k,m$ there is, up to scaling, a unique 
section $\tau_{k,m}$ of $\caL_{k,m}\otimes\caL_{k+1,m}^{-1}$ with
divisor $S-\tilde O_{k,m}$.  Using (\ref{eq:tmap})
this fixes an isomorphism
$t_{k,m}:\Gamma(\caL_{k+1,m})\to\Gamma(\caL_{k,m})$ which identifies sections by
their behaviour over $\lambda=\infty$. Similarly using (\ref{eq:tmap}) we fix an isomorphism $\hat
t_{k,m}:\Gamma(\caL_{k,m+1})\to\Gamma(\caL_{k,m})$ by choosing a section  
$\hat\tau_{k,m} $ of $\caL_{k,m}\otimes\caL_{k,m+1}^{-1}$ with divisor $Q-\tilde O_{k,m}$. 
Notice that $\hat t_{k,m}\circ t_{k,m+1}$ and $t_{k,m}\circ\hat t_{k+1,m}$ differ only by a scaling,
since $\hat\tau_{k,m}\otimes\tau_{k,m+1}$ and $\tau_{k,m}\otimes\hat\tau_{k+1,m}$ differ only by a
scaling (they are sections of the same line bundle and have the same divisor).
Therefore we obtain, by composition, a well-defined
isomorphism $\P\Gamma(\caL_{k,m})\simeq\P\Gamma(\caL)$. 
Since the only freedom in the choice of $\tau_{k,m}$ and $\hat\tau_{k,m}$ is the scale, 
which is irrelvant when we pass to
the projective spaces, this construction depends only on the spectral data.
Finally, by fixing an identification of $\P\Gamma(\caL)^*$ with $\P^1$ we obtain the maps
$\zeta_{k,m}$.  Since the last step is not determined by the spectral data the 
discrete map is only determined up to M\" obius equivalence.

At this point we can see why, under the conditions of the theorem, the map $z$ should be periodic. It
suffices to show that $\zeta_{k+n,m}=\zeta_{k,m}$. Consider the map 
$r_{k,m}=t_{k+n-1,m}\circ\ldots\circ t_{k,m}$: by (\ref{eq:tmap}) it is characterised by
\[
r_{k,m}:\Gamma(\caL_{k,m})\to\Gamma(\caL_{k,m}) \quad \hbox{where}\quad r_{k,m}(\sigma)\vert\infty = 
(\sigma\otimes h_{k,m})\vert\infty,
\]
for $h_{k,m}=\tau_{k+n-1,m}\otimes\ldots\otimes\tau_{k,m}$. Notice that $h_{k,m}$ is a rational
section of $\caL_{k,m}\otimes\caL_{k+n,m}^{-1}$ with divisor $\frac{n}{2}(2S-O-\tilde O)$ for $n$ even
and $n(S-O)$ for $n$ odd. Therefore, under the conditions of the theorem,
$h_{k,m}$ is a rational function on $\Sigma$ with $h_{k,m}(P_\infty)=h_{k,m}(\tilde P_\infty)$.
It follows that $r_{k,m}$ is a scalar multiple of the identity whence
$\zeta_{k+n,m}=\zeta_{k,m}$.

Now we must show that this recipe produces discrete conformal maps.
First we observe the following convenient result.
\begin{lem}
\label{lem:SQ}
The map $\zeta_{k,m}$ sends $\Gamma(\caL_{k,m}(-S))$ and $\Gamma(\caL_{k,m}(-Q))$ to
$z_{k+1,m}$ and $z_{k,m+1}$ respectively.
\end{lem} 
{\bf Proof.} It suffices to show that the maps $t_{k,m}$ and $\hat t_{k,m}$ in (\ref{eq:tmap}) send,
respectively, $\Gamma(\caL_{k+1,m}(-O_{k+1,m}))$ to $\Gamma(\caL_{k,m}(-S))$ and
$\Gamma(\caL_{k,m+1}(-O_{k,m+1}))$ to $\Gamma(\caL_{k,m}(-Q))$. We will prove the former: the latter use
the same proof with $S$ replaced by $Q$. Let $\sigma$ generate the line
$\Gamma(\caL_{k+1,m}(-O_{k+1,m}))$, then $\sigma$ is a holomorphic section with a zero at
$O_{k+1,m}=\tilde O_{k,m}$ and this is where $\tau_{k,m}$ has a simple pole. Hence 
$\sigma\otimes\tau_{k,m}$ is a
holomorphic section of $\caL_{k,m}$ with a zero at $S$ i.e.\ it generates $\Gamma(\caL_{k,m}(-S))$.
$\Box$

To prove the theorem we will show
that the action of tensoring sections with $\tau_{k,m}$ is represented 
by a matrix of the form (\ref{eq:Tmatrix}) and
deduce from this that the cross-ratios are constant. Since this matrix really acts on the
eigenline bundle $\caE_{k,m}$ dual to $\caL_{k,m}$ we must take some care to 
describe the relationship
between rational sections of $\caE_{k,m}$ and global sections of $\caL_{k,m}$. 
\begin{lem}
\label{lem:dual}
Let $\caL$ be a line bundle of degree $g+1$ over the genus $g$ hyperelliptic curve $\Sigma$ 
with $\Gamma(\caL(-P-\tilde P))=0$ (for some $P\in
\Sigma$) and let $\caE$ denote the dual line bundle. Then $\Gamma(\caL)$ is 
canonically dual to $\Gamma(\caE(R))$ where $R$ is the ramification divisor of
$\lambda$. This duality identifies the hyperplane $\Gamma(\caL(-P))$ with the line
$\Gamma(\caE(R-\tilde P))$.
\end{lem}
{\bf Proof.} Let $\caF$ be the field of rational functions on $X$ and let $\caK$ denote the subfield of
rational functions of $\lambda$.
Clearly $\caF$ is a two dimensional $\caK$-space and we have a $\caK$-linear
map $\Tr:\caF\to\caK$ which gives to each element $f\in\caF$ the trace of the matrix in $\mathfrak{
gl}_2(\caK)$ representing
the multiplication map $a\mapsto fa$ on $\caF$ (since the trace is invariant this is independent of
the $\caK$-basis chosen for $\caF$). It is easy to check that:
(a) $\Tr(f)$ is globally holomorphic
precisely when its divisor of poles is no worse than the ramification divisor $R$ of $\lambda$; 
(b) for any such function $f$, $\Tr(f)=0$ if its divisor of zeroes includes $P+\tilde P$ (for some
$P\in\Sigma$).
We have a non-degenerate
$\caK$-bilinear form on $\caF(\caE)\times\caF(\caL)$ by $(v,e)\mapsto \Tr(e(v))$ which, by 
the properties (a) and (b), pairs $\Gamma(\caL)$ non-degenerately with
$\Gamma(\caE(R))$. Further, if $e\in\Gamma(\caL(-P))$ then $\Tr(e(v))=0$ if and only if 
$v\in\Gamma(\caE(R-\tilde P))$.$\Box$

Since $\caE(R)$ also has the property $\Gamma(\caE(R-P-\tilde P))=0$ it follows 
that $\caF(\caE))=\caK\otimes_\C\Gamma
(\caE(R))$. As a result we may draw the following commuting diagram with which we define the  map $T_{k,m}$:
\begin{equation}
\label{eq:Tdiagram}
\begin{array}{ccc}
\caF(\caE_{k,m}) & \stackrel{\otimes\tau_{k,m}}{\rightarrow} & \caF(\caE_{k+1,m})\\
\downarrow       &           & \downarrow \\
\caK\otimes\Gamma(\caL_{k,m})^* &\stackrel{T_{k,m}}{\rightarrow} &\caK\otimes\Gamma(\caL_{k+1,m})^*
\end{array}
\end{equation}
Now when we identify $\Gamma(\caL_{k,m})^*$ with $\C^2$ using any lift of $\zeta_{k,m}$
we obtain a
$\caK$-valued matrix we will define to be $T^\lambda_{k,m}$. The theorem will follow from the next
proposition.
\begin{prop}
\label{pp:T}
The matrix $T^\lambda_{k,m}$ obtained from (\ref{eq:Tdiagram}) is of the form $I-\lambda^{-1}A_{k,m}$
where: (i) $\ker(A_{k,m})=z_{k,m}$; (ii) $\im(A_{k,m})=z_{k+1,m}= \ker(I-A_{k,m})$; and (iii)
$I-q^{-1}A_{k,m}$ maps $z_{k,m+1}$ to $z_{k+1,m+1}$ for $q=\lambda(Q)$.   
It follows that $[z_{k,m+1}:z_{k,m}:z_{k+1,m}:z_{k+1,m+1}] = q$ and therefore the
map $z$ is discrete conformal.
\end{prop}
{\bf Proof.} To begin, observe that $T^\infty_{k,m}$ represents the map $\sigma\to\sigma\otimes\tau_{k,m}$
over $\lambda=\infty$. By (\ref{eq:tmap}) this is the identity. Further, since $\tau_{k,m}$ has degree
one it follows that there is a constant matrix $A_{k,m}$ such that $T^{\lambda}_{k,m} = I-\lambda^{-1}
A_{k,m}$. To prove (i), (ii) and (iii) we consider three different values of $\lambda$. 

First, $A_{k,m}
= (\lambda T^\lambda_{k,m})\vert_{\lambda=0}$ represents $\lambda\tau_{k,m}$ over $\lambda=0$. 
But $\lambda\tau_{k,m}$ has
a simple zero at $O_{k,m}$ and none at $\tilde O_{k,m}$, hence $(\sigma\otimes\lambda\tau_{k,m})\vert
_{\lambda=0}=0$ if and only if $\sigma$ vanishes at $\tilde O_{k,m}$. By lemma \ref{lem:dual}
$z_{k,m}$ corresponds to $\Gamma(\caE_{k,m}(R-\tilde O_{k,m}))$ and $z_{k+1,m}$
corresponds to $\Gamma(\caE_{k+1,m}(R- O_{k,m}))$ (since $\tilde O_{k+1,m}=O_{k,m}$)
so $\ker(A_{k,m})=z_{k,m}$ while $\im(A_{k,m}) = z_{k+1,m}$.

Secondly, consider $I-A_{k,m}$, which represents $\tau_{k,m}$ over $\lambda=1$. But $\tau_{k,m}$
has a zero at $S$ so
$(\sigma\otimes\tau_{k,m})\vert_{\lambda =1}=0$ if and only if $\sigma$ has a zero at $\tilde S$.Recall
that 
$z_{k+1,m}$ corresponds to $\Gamma(\caE_{k,m}(R- \tilde S))$ by lemmas \ref{lem:SQ} and \ref{lem:dual}
so $z_{k+1,m} = \ker(I-A_{k,m})$.

Thirdly, $I-q^{-1}A_{k,m}$ represents $\tau_{k,m}$ over $\lambda=q$. Here $\tau_{k,m}$ has neither 
zeroes nor
poles so $\sigma\otimes\tau_{k,m}$ has a zero at $\tilde Q$ if and only if $\sigma$ does. Since
$z_{k,m+1}$ corresponds to $\Gamma(\caE_{k,m}(R- \tilde Q))$ and
$z_{k+1,m+1}$ corresponds to $\Gamma(\caE_{k+1,m}(R- \tilde Q))$ 
(by lemmas \ref{lem:SQ} and \ref{lem:dual}) we see that $I-q^{-1}A_{k,m}$ maps $z_{k,m+1}$ to 
$z_{k+1,m+1}$.

Finally, it is an elementary computation (which we will leave to the reader) to establish that for
$A_{k,m}$ to have 
these properties we are obliged to have $[z_{k,m+1}:z_{k,m}:z_{k+1,m}:z_{k+1,m+1}] = q$.$\Box$

\section{The Lax pair and a loop group action.}

In \cite{NijC} Nijhoff and Capel wrote down a Lax pair for the discrete conformal 
map equations. Using this we will
show that all discrete conformal maps: (i) come in 1-parameter families obtained by deforming the
cross-ratio (we have already seen this for periodic maps in corollary \ref{cor:family});
(ii) admit non-trivial transformations by
elements of an infinite dimensional Lie group (in fact, a loop group). To begin, let us define
$\caZ_q=\{z:\Z^2\to\C\vert$ discrete conformal with cross-ratio $q$ and $z_{0,0}=0\}$. We will say 
maps of this type are based at 0.  Since we are working with maps into the plane we may define
\begin{equation}
\label{eq:uv}
u_{k,m}= z_{k+1,m}-z_{k,m}\ , \qquad
v_{k,m}=z_{k,m+1} - z_{k,m}.
\end{equation}
Now fix a pair $\alpha,\beta\in \C$ for which $\beta^2/\alpha^2=q$ and define
\begin{equation}
\label{eq:UV}
U_{k,m}(\lambda) = \left(\begin{array}{cc} 1 & \alpha^2 u_{k,m}^{-1}\lambda^{-1}\\
                                     u_{k,m} & 1\end{array}\right)\ ,
\qquad
V_{k,m}(\lambda) = \left(\begin{array}{cc} 1 & \beta^2 v_{k,m}^{-1}\lambda^{-1}\\
                                     v_{k,m} & 1\end{array}\right).
\end{equation}
One readily computes that the pair $U,V$ is a discrete Lax pair i.e.
\begin{equation}
\label{eq:Lax}
U_{k,m}V_{k+1,m} = V_{k,m}U_{k,m+1}.
\end{equation}
We consider these as maps $\Z^2\to\caG$ where $\caG$ is the infinite dimensional 
Lie group of holomorphic maps from
$\P^1_\lambda\backslash\{0,\alpha^2,\beta^2\}$ to $GL_2$. By the Lax equations 
(\ref{eq:Lax}) there exists a unique map $\Phi:\Z^2\to\caG$ satisfying:
\begin{equation}
\label{eq:Phi}
\Phi_{k+1,m} = \Phi_{k,m}U_{k,m}\ ; \qquad
\Phi_{k,m+1} = \Phi_{k,m}V_{k,m}\ ; \qquad \Phi_{0,0}=I.
\end{equation}
We will call such maps {\em extended frames}.  Each $z\in\caZ_q$ has a unique extended frame $\Phi$ and
$z$ is recovered from the first column of
$\Phi(\infty)$. To see this simply observe that
$\Phi$ is constructed inductively from (\ref{eq:Phi}) and is uniquely
specified by the initial condition $\Phi_{0,0}=I$. One readily checks that
\[
\left(\begin{array}{cc} 1 & 0\\
                  z_{k,m} & 1
\end{array}\right)
\]
satisfies (\ref{eq:Phi}) for $\lambda=\infty$ so, by uniqueness, this must be $\Phi_{k,m}(\infty)$.

\noindent
{\bf Example: the vacuum solution.} The simplest example of a discrete conformal map is the tiling of the
plane by a parallelogram of cross-ratio $q$ i.e.\
\[
z_{k,m}=k\alpha +m\beta .
\]
We will adopt a common terminology from soliton theory and call this solution the
`vacuum solution': it corresponds to having $u_{k,m}=\alpha,v_{k,m}=\beta$ for all $k,m$.
In this case it is easy to check that for this map the extended frame is
\begin{equation}
\label{eq:vacuum}
\Phi_{k,m} = (I+\alpha \Lambda)^k(I+\beta \Lambda)^m\ , \quad \hbox{where}\ 
\Lambda=\left(\begin{array}{cc}
0 & \lambda^{-1}\\ 1 & 0\end{array}\right).
\end{equation}

\smallskip
Now let us show that each $z\in\caZ_q$ gives rise to a 1-parameter family of discrete conformal maps. 
Set $e_1=(1,0)^t, e_2=(0,1)^t$ and let $B_0\subset GL_2$ denote the subgroup of upper triangular matrices: this is the
stabilizer of the line generated by $e_1$. For each 
$\lambda\in\P^1_\lambda\backslash\{0,\alpha^2,\beta^2\}$ define
$z(\lambda):\Z^2\to\P^1$ by taking $z_{k,m}(\lambda)$ to be the line generated by 
$F_{k,m}(\lambda)=\Phi_{k,m}(\lambda)e_1$. We think of $F(\lambda):\Z^2\to\C^2$ as a lift of $z(\lambda)$.
\begin{prop}
The map $z(\lambda)$ is discrete conformal with cross-ratio
$\beta^2(1-\lambda^{-1}\alpha^2)/\alpha^2(1-\lambda^{-1}\beta^2)$. If $z,\hat z\in \caZ_q$ 
are M\" obius equivalent
maps then so are $z(\lambda),\hat z(\lambda)$ for each $\lambda$.
\end{prop}
{\bf Proof.} The cross-ratio of the map $z(\lambda)$ may be written using the lift $F(\lambda)$ as
\[
\frac{\det(F_{k,m+1}(\lambda)\ F_{k,m}(\lambda))}{\det(F_{k,m}(\lambda)\
F_{k+1,m}(\lambda))}
\frac{\det(F_{k+1,m}(\lambda)\ F_{k+1,m+1}(\lambda))}{\det(F_{k+1,m+1}(\lambda)\
F_{k,m+1}(\lambda))}.
\]
where if $a,b\in\C^2$ then $(a\ b)$ denotes the matrix with these columns.
But $F_{k+1,m} = \Phi_{k,m}U_{k,m}e_1$ and $F_{k,m+1} =
\Phi_{k,m}V_{k,m}e_1$. Inserting these expressions and 
using (\ref{eq:UV}), (\ref{eq:Lax}) this reduces to 
\[
\frac{
\det\left(\begin{array}{cc} 1 & 1 \\ v_{k,m} & 0\end{array}\right)}
{\det\left(\begin{array}{cc} 1 & 1 \\ 0 & u_{k,m} \end{array}\right)}
\frac{
\det\left(\begin{array}{cc} 1 & 1 \\ 0 & v_{k+1,m} \end{array}\right)}
{\det\left(\begin{array}{cc} 1 & 1 \\ u_{k,m+1} & 0\end{array}\right)}
\frac{\det U_{k,m}}{\det V_{k,m}}
=\frac{\beta^2(1-\lambda^{-1}\alpha^2)}{\alpha^2(1-\lambda^{-1}\beta^2)},\ 
\]
Now suppose $z, \hat z$ are M\" obius equivalent. Since they are both based at $0$ it suffices to consider 
two types of M\" obius transforms: 1) $\hat z =cz$; 2) $\hat z = z/(1-cz)$ for any non-zero constant $c$. 
We wish to show that in each case there exists $a:\Z^2\to B_0$
for which 
\[
\hat U_{k,m} = a^{-1}_{k,m}U_{k,m}a_{k+1,m},\quad \hat V_{k,m}=a_{k,m}^{-1}V_{k,m}a_{k+1,m}.
\]
For then the extended frames $\Phi,\hat\Phi$ satisfy $\hat\Phi_{k,m}=
a_{0,0}^{-1}\Phi_{k,m}a_{k,m}$, whence $\hat
z(\lambda)=a_{0,0}^{-1}\circ z(\lambda)$ (M\" obius action). 
We leave it to the reader to verify that for the two cases above we can take:
\[
1)\ a_{k,m} = \left(\begin{array}{cc} c& 0\\ 0& 1\end{array}\right), \qquad
2)\ a_{k,m} = \left(\begin{array}{cc} z_{k,m}-c^{-1}& -1\\ 0& 1/(z_{k,m}-c^{-1})\end{array}\right).
\Box
\]
{\bf Remark.} Notice that each $z(\lambda)$ will be based at 0, although we cannot guarantee that it maps into
$\C\subset\P^1$. The parameterisation of this family $z(\lambda)$ depends upon the choice of $\alpha^2,\beta^2$
but the family itself clearly does not, since a rescaling of $\lambda$ will allow us to move through all pairs
$\alpha,\beta$ such that $\beta^2/\alpha^2=q$. 

\smallskip
\begin{cor}
Any Lax pair of the form (\ref{eq:UV}) produces a 1-parameter family of discrete conformal maps.
\end{cor}
This follows immediately from the first part of the previous proof.

\subsection{The dressing action.}

Let us now fix $q$ and choose $\alpha,\beta$ so that $\vert\alpha\vert,\vert\beta\vert < 1$ 
(such a pair can clearly be found for each $q$). Then we may
view $U,V$ (and hence $\Phi$) as taking values in the connected loop group $LGL_2 = 
\{g\in C^\omega(S^1,GL_2):\det(g)$ has winding number zero$\}$. In fact 
$U,V$ and $\Phi$ all take values in the subgroup $N=\{g\in LGL_2:$ $g$ extends
holomorphically into the disc $\vert\lambda^{-1}\vert <1$ with $g(\infty)$ lower unipotent $\}$ (we will say 
a matrix $A$ is lower unipotent if $I-A$ is strictly lower triangular). 
The dressing action uses the fact that almost
every $g\in LGL_2$ can be factorised into a product $g=g_N.g_B$ where $g_N\in N$ and $g_B$ belongs to the
subgroup $B=\{g\in LGL_2:$ $g$ extends holomorphically 
into $\vert\lambda\vert<1$ with $g(0)\in B_0\}$. One knows (from e.g.\
\cite{PreS}) that the space $N.B$ of all products is open dense in
$LGL_2$. Given an extended frame $\Phi$ and $g\in B$ 
we might expect $g\Phi$ to map into $N.B$, in which case we could define 
$g\circ\Phi=(g\Phi)_N$ and therefore we have
\begin{equation}
\label{eq:factor}
g\Phi_{k,m}=(g\circ\Phi)_{k,m}\Psi_{k,m}\ ,\quad\hbox{for some}\ \Psi:\Z^2\to B.
\end{equation} 
\begin{lem}
For any $g\in B$ for which it exists, the map $g\circ\Phi:\Z^2\to N$ is again the extended frame for a discrete
conformal map $g\circ z$ of cross-ratio $\beta^2/\alpha^2$.
\end{lem}
{\bf Proof.} It suffices to show that the Lax pair for $g\circ\Phi$ is again of the form
(\ref{eq:UV}), since it is clear that $g\circ\Phi_{0,0}=g_N=I$. Observe that 
$g\circ\Phi=g\Phi\Psi^{-1}$ so that we can define
\[
\begin{array}{c}
\hat U_{k,m}= (g\circ\Phi)_{k,m}^{-1}(g\circ\Phi)_{k+1,m} = \Psi_{k,m}U_{k,m}\Psi_{k+1,m}^{-1}\ ,\\
\hat V_{k,m}= (g\circ\Phi)_{k,m}^{-1}(g\circ\Phi)_{k,m+1} = \Psi_{k,m}V_{k,m}\Psi_{k,m+1}^{-1}\ .
\end{array}
\]
Since $\Psi_{k,m}\in B$ we can expand it in Fourier series as
\[
\Psi_{k,m}=\left(\begin{array}{cc} a_{k,m} & b_{k,m}\\ 0 & c_{k,m}\end{array}\right) + O(\lambda).
\]
Consequently we compute
\[
\Psi_{k,m}U_{k,m}\Psi_{k+1,m}^{-1} = \alpha^2\lambda^{-1}\left(\begin{array}{cc} 0& e_{k,m}\\ 0&0\end{array}\right)
+O(1),
\]
for some expression $e_{k,m}$. But $g\circ\Phi$ takes values in $N$ so the only other terms in
this expression are constant in $\lambda$ and lower unipotent i.e.
\[
\hat U_{k,m} = \left(\begin{array}{cc} 1 & \alpha^2\lambda^{-1}e_{k,m} \\ \hat u_{k,m} & 1\end{array}\right),
\]
for some $\hat u_{k,m}$. Now we claim that $\det(\Psi_{k,m})=\det(g)$ for all $k,m$, whence $\det(\hat U_{k,m}) =
\det(U_{k,m})=1-\alpha^2\lambda^{-1}$ so that $e_{k,m}=\hat u_{k,m}^{-1}$ and we are done. To see the claim,
notice that
\[
\det(\Phi)\det(g\circ\Phi)^{-1}=\det(g)^{-1}\det(\Psi).
\]
The terms on the left extend holomophically into $\vert\lambda\vert <1$ while on the right they extend 
holomorphically into
$\vert\lambda^{-1}\vert <1$. Therefore both sides are constant and by evaluating at $k,m=0$ we see this
constant is 1. A similar argument for $\hat V$ finishes the proof.$\Box$

Formally, at least, one can check that $\Phi\mapsto g\circ\Phi$ is a (left) group action and therefore we
obtain an action of $B$ on $\caZ_q$.  By analogy with the theory of the KdV equation we call this 
the dressing action (see e.g.\cite{Wil}). The obstruction to making this more rigorous 
is that we cannot guarantee that $g\circ\Phi$ exists off $(0,0)\in\Z^2$ (compare
with the KdV theory where one has at least local solutions guaranteed). 

\noindent
{\bf Remark: the effect of the scaling $\alpha,\beta\mapsto k\alpha,k\beta$.}
For fixed $q$ we are only permitted to change $\alpha,\beta$ by the scaling 
$k\alpha,k\beta$, where we require 
$\vert k\vert <\min\{\vert\alpha\vert^{-1},\vert\beta\vert^{-1}\}$. 
If $U(\lambda),V(\lambda)$ is the Lax pair for $z,\alpha,\beta$ then
the Lax pair for $z,k\alpha,k\beta$ is $U(k^{-2}\lambda),V(k^{-2}\lambda)$. 
It is clear that our choice of working
with loops on the unit circle limits us if the scaling is not unimodular. But this limitation is unnecessary: we
could equally well work with loops on any circle and apply the dressing theory. It is not hard to see that by
rescaling the unit circle
the dressed map $g\circ z$ obtained using $\alpha,\beta$ can be obtained using $k\alpha,k\beta$ and
$g(k^{-2}\lambda)$.  
 
\smallskip
For the next result, we think of $B_0$ as the subgroup of constant loops in $B$.
\begin{lem} 
The dressing action of the subgroup $B_0\subset B$ corresponds precisely to the action of $B_0$ as the full
group of base point preserving M\" obius transformations.
\end{lem}
{\bf Proof.} Let $z\in\caZ_q$ and let $\Phi$ be its extended frame, so $F=\Phi(\infty)e_1$ lifts $z$ into
$\C^2$. Take $g\in B_0$, then $g\Phi$ extends holomorphically into $\vert\lambda^{-1}\vert <1$. 
If $g\Phi_{k,m}$ factorises then $g\circ
z$ (i.e.\ the action of $g$ by M\" obius transforms) has lift
\[
gF = (g\Phi)_N(\infty)(g\Phi)_B(\infty)e_1 =k(g\circ\Phi)(\infty)e_1, \quad\hbox{for some $k\in\C$,}
\]
since $(g\Phi)_B(\infty)$ takes value in $B_0$ which is the stabilizer of the line generated by
$e_1$. Now we observe that $g\Phi_{k,m}$ factorises precisely when $(g\Phi_{k,m})(\infty)e_1\neq e_2$, i.e.\ when 
$z_{k,m}$ belongs to the plane $\C$.
$\Box$

Now observe that $B_0$ is a normal subgroup of $B$ and $B/B_0\simeq LG_+ = \{g\in B: g(0) = I\}$. Since the
quotient space $\caZ_q/B_0$ can be identified with the space of all M\" obius equivalence classes of
discrete conformal maps we have an induced action of $LG_+$ on $\caZ_q/B_0$. This is what we 
shall focus on from now on. 

\subsection{Dressing orbit of the vacuum solution.}

Our main interest in the dressing action is to examine the dressing orbit of the vacuum solution,
which we will denote by $z^{(0)}$.
Throughout this section $\Phi$ will be the extended frame (\ref{eq:vacuum}) for that solution.
First we will describe this dressing orbit as a quotient space. It can be shown (by a tedious but
straightforward argument using Fourier expansions) that $g\circ\Phi=\Phi$ if and only if $g\in\Gamma_+=\{g\in
B:g\Lambda=\Lambda g\}$. Therefore the dressing orbit of $\Phi$ is identifiable with $B/\Gamma_+$, so the
dressing orbit of the vacuum solution in $\caZ_q/B_0$ is identifiable with the double coset space
$B_0\backslash B/\Gamma_+$. 
Since $B_0$ is a normal subgroup the group $LG_+$ acts transitively on the left so we may identify 
$B_0\backslash B/\Gamma_+$ with
the quotient space $LG_+/\Gamma_+$. Here the action of $\Gamma_+$ on $LG_+$ is given by $\gamma\cdot 
g=\gamma(0)^{-1} g\gamma$.

The principal result in this section is that this dressing orbit contains all those
periodic discrete conformal maps which are "small perturbations" of their corresponding SKdV solution
in the following sense.
\begin{thm}
\label{th:orbit}
Let $z:\Z^2\to \P^1$ be a periodic discrete conformal map (of cross-ratio $q$)
with spectral data $(\Sigma,\lambda,O,S,Q)$ for
which: (a) $\lambda =0$ is a branch point (i.e.\ $O$=$\tilde O$); (b) the disc $\vert\lambda\vert 
\leq\max\{ 1,\vert q\vert\}$ contains no other branch points. Then $z$ is in the dressing orbit of the vacuum
solution.
\end{thm}
To prove this we will use the Grassmannian picture of the dressing orbit of the KdV vacuum solution
developed by Segal and Wilson \cite{SegW}(see also \cite{PreS}). 
Let $H=L^2(S^1,\C)$ and recall that $LGL_2$ acts on this space
by `interleaving Fourier series' i.e.\ we use the isometric isomorphism between $H$ and $L^2(S^1,\C^2)$
given by $ f(\zeta) \mapsto (f_0(\lambda),f_1(\lambda))^t$ where $f_0(\zeta^2) +
\zeta^{-1}f_1(\zeta^2)=f(\zeta)$.
The usual action of $LGL_2$ on $\C^2$-valued functions then passes to an action 
on $H$ which we will denote by $g\cdot f$. In particular,
the action of the commutative subgroup $\Gamma=\{ g\in LGL_2:g\Lambda=\Lambda g\}$ on $H$ corresponds to
the multiplication action of $C^\omega(S^1,\C^*)$. This follows from the observation that
$\Lambda\cdot f(\zeta) = \zeta^{-1}f(\zeta)$. Therefore we have
\begin{equation}
\label{eq:scalar}
\Phi_{k,m}\cdot f(\zeta) = (1+\alpha\zeta^{-1})^k(1+\beta\zeta^{-1})^mf(\zeta)
\end{equation}
for the vacuum extended frame.

Next, recall that $H$ has an orthogonal decomposition $H_-\oplus H_+$ where $H_-$ (respectively, $H_+$)
is the subspace of functions whose Fourier series possess only non-positive (respectively, positive)
powers of $\zeta$. [Readers familiar with \cite{PreS,SegW} should note that we have applied the involution
$\lambda\mapsto\lambda^{-1}$ to the the picture described there.] We will define the Grassmannian to be
$Gr = \{ g\cdot H_-:g\in LGL_2\}$. We recall that the orbit of $H_-$ under the subgroup $LG_+$ is open
dense (the `big cell') and that it is characterised as being the set of all $W\in Gr$ for which the
orthogonal projection $pr_-:W\to H_-$ is invertible. Our assumption about the existence of $g\circ\Phi$
implies that $\Phi_{k,m}^{-1}\cdot W$ belongs to the big cell for all $k,m$, where
$W=g^{-1}\cdot H_-$.

We will now associate to each $W\in Gr$ the discrete analogue of the Baker function used in
\cite{SegW}. First we observe that we may
refine our current factorisation $N.B$ into $N.T.\tilde N$ where $T$ denotes the subgroup of constant
diagonal loops and $\tilde N=\{g\in B:g(0)$ is upper unipotent$\}$. Thus any $g\in N.B$ can be factorised
into $g_Ng_Tg_{\tilde N}$.  Let us now fix a $g\in LG_+$ and define 
\[
\psi = \Phi\Psi_{\tilde N}^{-1}\cdot 1 = g^{-1}(g\circ\Phi)\Psi_T\cdot 1,
\]
using (\ref{eq:factor}). Since both subgroups $N$ and $T$ preserve $H_-$ we see that 
$\psi:\Z^2\to W=g^{-1}\cdot H_-$.
\begin{lem}
\label{lem:WBaker}
Each $\psi_{k,m}\in W$ given in this way is uniquely determined by the property that
\[
pr_-((1+\alpha\zeta^{-1})^{-k}(1+\beta\zeta^{-1})^{-m}\psi_{k,m})=1.
\]
\end{lem}
{\bf Proof.} Since we are assuming $\Phi^{-1}\cdot W$ is always in the big cell it suffices to show that
$pr_-(\Phi_{k,m}^{-1}\cdot\psi_{k,m})=1$, given (\ref{eq:scalar}). But $\Phi^{-1}\cdot\psi = 
\Psi_{\tilde N}^{-1}\cdot 1$. It is easy to check that for any $n\in \tilde N$, 
$n\cdot 1 = 1+ O(\zeta)$.$\Box$

We will call this function the discrete Baker function for $W$.
Now we want to reconstruct the dressed map $g\circ z^{(0)}$ from $\psi$. First we
observe that $\zeta^{-2}W\subset W$  and the quotient $W/\zeta^{-2}W$ is two dimensional (since
$g\cdot\zeta^{-2}f=\zeta^{-2}g\cdot f$ and $H_-/\zeta^{-2}H_-$ is two dimensional).
\begin{lem}
\label{lem:quotient}
Choose any linear identification of $W/\zeta^{-2}W$ with $\C^2$ and let $[\psi]:\Z^2\to\P^1$ be the
map given by $(k,m)\mapsto \psi_{k,m} +\zeta^{-2}W$. Then 
$[\psi]$ and $g\circ z^{(0)}$ are M\" obius equivalent
(so in particular this gives a discrete conformal map).
\end{lem}
{\bf Proof.} Observe that $1+\zeta^{-2}H_-$ and $\zeta^{-1}+\zeta^{-2}H_-$ span $H_-/\zeta^{-2}H_-$.
Therefore by setting $e_1=g^{-1}\cdot 1$ and $e_2=g^{-1}\cdot \zeta^{-1}$ we obtain a basis
$e_1+\zeta^{-2}W,e_2+\zeta^{-2}W$ for $W/\zeta^{-2}W$. We can write $g\circ\Phi$ and $\Psi_T$
in the form
\[
g\circ\Phi = \left(\begin{array}{cc} a & 0 \\ b & c \end{array}\right) + O(\lambda^{-1}),\quad
\Psi_T=\left(\begin{array}{cc} s & 0 \\ 0 & *\end{array}\right)
\]
so that the map $g\circ z^{(0)}$ has homogeneous coordinates $[a,b]$ and $\psi_T\cdot 1=s$. 
Now observe that 
$(g\circ\Phi)\Psi_T\cdot 1 = sa +\zeta^{-1} sb + O(\zeta^{-2})$ 
and therefore 
\[
\begin{array}{lll} s(ae_1 +be_2) & \equiv & g^{-1}(g\circ\Phi)\Psi_T\cdot 1 \bmod \zeta^{-2}W\\
                              & \equiv & \psi \bmod \zeta^{-2}W. \end{array}        
\]
It follows that in this basis the map $[\psi]$ has homogeneous coordinates $[a,b]$.$\Box$

Now we are in a position to prove theorem \ref{th:orbit}. Given a periodic discrete conformal map
$z$ with spectral data satisfying the conditions of the theorem we will construct a $W\in Gr$ for
which the map $\psi:\Z^2\to W$ recovers $z$ as $[\psi]$. By the previous lemma this will prove the
theorem. 

Given $(\Sigma,\lambda,\caL,O,S,Q)$ we follow \cite{SegW} and define $W$ in the following manner. By
assumption there exists $r > \max\{1,\vert q\vert\}$ for which the disc $\vert\lambda\vert\leq r$
contains no branch points. We set $\zeta=\sqrt{\lambda/r}$ so that the unit $\zeta$-circle may be
identified with the boundary of a disc $\Delta$ 
on $\Sigma$ about $O$ which contains no ramification points
other than $O$ i.e.\ this is a coordinate disc with coordinate $\zeta$ and contains the points
$O,S,Q$. We define $W$ to be the $L^2$-closure of the space $W^0$ of all $f\in H$ which extend
meromorphically into $\Sigma -\Delta$ where they have divisor of poles no worse than $D$, where
this is the unique positive divisor of degree $g$ in the class $\caL(-O)$. One knows that $W\in Gr$
(even the small Grassmannian we have chosen, since the whole construction is analytic). Now, with the
definitions from previous sections, we define $\psi_{k,m}$ to be (the boundary of) a meromorphic
function with divisor (\ref{eq:divisor}). Since the points $O,S,Q$ are not 
in $\Sigma-\Delta$ this function
belongs to $W$ for all $k,m$. Moreover, if we set $\alpha=-\zeta(S)$ and $\beta =-\zeta(Q)$ we see that
\[
(1+\alpha\zeta^{-1})^{-k}(1+\beta\zeta^{-1})^{-m}\psi_{k,m}
\]
extends holomorphically into $\Delta$ and $\psi_{k,m}$ will be uniquely determined by requiring this
expression to equal $1$ when evaluated at $\zeta=0$. By lemma \ref{lem:WBaker} this must be the discrete
Baker function for $W$ and by lemma \ref{lem:qBaker} the map
$(k,m)\to\psi_{k,m}(P_\infty)/\psi_{k,m}(\tilde P_\infty)$ recovers $z:\Z^2\to \P^1$ up to M\" obius
equivalence. But the map $W^0\to \P^1$ given by $f\mapsto [f(P_\infty),f(\tilde P_\infty)]$
induces an isomorphism $\P(W^0/\zeta^{-2}W^0)\to\P^1$. 
Since $W^0$ is dense in $W$ the spaces $W^0/\zeta^{-2}W^0$ and
$W/\zeta^{-2}W$ are equal. Therefore, by the previous lemma the maps $[\psi]$ and $z$ are equivalent,
hence $z$ is in the dressing orbit of the vacuum solution.

\smallskip\noindent
{\bf The discrete cubic.} 
Lemma \ref{lem:quotient} allows us to compute examples which are difficult to obtain
otherwise. For example, take $W=\C\langle\zeta\rangle +\zeta^{-1}H_-$, where $\C\langle\zeta\rangle$ is the
vector space generated by $\zeta$. One knows that $W\in Gr$
(see \cite{SegW} \S 7): indeed $W$ lies in the Grassmanian for rational loops in $GL_2$. It is
elementary to show that in this case
\[
\psi_{k,m} = (1+\alpha\zeta^{-1})^k(1+\beta\zeta^{-1})^m(1-\frac{\zeta}{k\alpha + m\beta}).
\]
Since $W/\zeta^{-2}W\simeq\C\langle\zeta,\zeta^{-2}\rangle$ we may take $\zeta + \zeta^{-2}W$ and
$\zeta^{-2}+\zeta^{-2}W$ as a basis. A straightforward computation shows that
\[
\psi_{k,m}\equiv \frac{1}{k\alpha+m\beta}\zeta - 
 \frac{z_{k,m}} {6(k\alpha +m\beta)}\zeta^{-2}\bmod \zeta^{-2}W,
\]
where
\[
z_{k,m} = 
    (k+1)k(k-1)\alpha^3 +3k^2m\alpha^2\beta +3km^2\alpha\beta^2 + (m+1)m(m-1)\beta^3.
\]
It follows that this is a discrete conformal map with cross-ratio $\beta^2/\alpha^2$. 
Indeed this is the discrete analogue of the cubic $z(x)=x^3$. In
the smooth (binomial) limit $(1+\alpha\zeta^{-1})^k(1+\beta\zeta^{-1})^m\to
\exp(x\zeta^{-1})$ we obtain the smooth Baker function $\psi(x)=\exp(x\zeta^{-1})(1-\zeta/x)$ and 
\[
\psi(x)\equiv -\frac{1}{x}\zeta + \frac{2}{3}x^2\zeta^{-2}\bmod \zeta^{-2}W.
\]
See \cite{SegW} \S 7 for a description of how this corresponds to taking the spectral curve
$\mu^2=\lambda^{-3}$ and for the similar spaces corresponding to the spectral curves with equations of
the form $\mu^2=\lambda^{-2g-1}$. All of these should produce (discrete) rational functions but it is
difficult to describe a priori which rational functions arise in this way.

\begin{figure}[ht]
\centering
\includegraphics[scale=0.5]{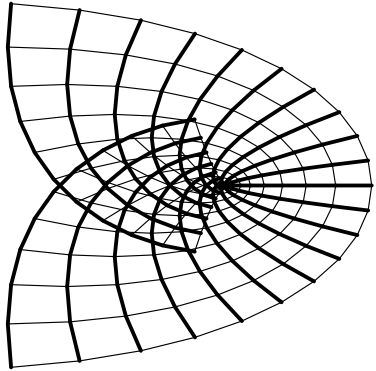}
\caption{The discrete cubic.}
\end{figure}

\subsection{Discrete conformal maps as Darboux transformations of the KdV hierarchy.}

In this section we will prove rigourously a rather surprising relationship between discrete conformal
maps and Darboux (B\" acklund) transforms of solutions to the KdV equation. We will work inside the
dressing orbit of the vacuum solution of KdV. Recall from \cite{SegW} that to every $W\in Gr$ there 
exists a solution $u_W(x,t)$ of the KdV equation, given by
\[
(\partial^2 +u_W)\psi_W=\zeta^{-2}\psi_W\ , \qquad \partial = \partial/\partial x,
\]
where $\psi_W:\C^2\to W$ is the (smooth) Baker function for $W$ i.e.\ the unique function with the property
that $pr_-(\exp(-x\zeta^{-1}-t\zeta^{-3})\psi_W)=1$ for almost all $x,t\in\C$. If $\Phi_{k,m}$ denotes
the vacuum frame (\ref{eq:vacuum}) then $W_{k,m}=\Phi_{k,m}^{-1}\cdot W$ defines a map $\Z^2\to Gr$.
We will show that $\{u_{W_{k,m}}\}$ is a $\Z^2$-family of Darboux transforms 
of $u_W$. 

To begin, let us recall, from e.g.\ \cite{EhlK} the Darboux transform for the KdV
hierarchy. Given any solution
$u(x,t)$ of the KdV equation (or any other equation in the hierarchy) we can produce other solutions
using the following procedure. Set $L=\partial^2+u$. 
For any pair $(\psi,c)$, consisting of a function $\psi(x,t)$ and a complex number $c$ 
for which $L\psi = c\psi$, 
define $v=\psi_x\psi^{-1}$. One readily checks that $u=-v_x-v^2+c$ and therefore
\[
L-c=(\partial +v)(\partial -v).
\]
The Darboux transform for the pair $(\psi,c)$ maps $u$ to $\tilde u$ where
\[
\tilde L-c=(\partial -v)(\partial +v)\ , \quad \hbox{i.e.}\ \tilde u = v_x-v^2+c.
\]
One can think of this as a formal conjugation of $L-c$ into $\tilde L-c$ in the algebra of
pseudo-differential operators and deduce that $\tilde u$ too satisfies the KdV equation (see
\cite{EhlK}). 
\begin{prop}
Whenever $0<\vert \alpha\vert <1$ the map $W\mapsto V=(1 +\alpha\zeta^{-1})W$ on $Gr$ induces the
Darboux transform on $u_W$ determined by taking $\psi=\psi_W(x,t;-\alpha)$ and $c=\alpha^{-2}$.
\end{prop}
{\bf Proof.} First, since $W$ is a linear space we can write $V=(\zeta^{-1}+\alpha^{-1})W$. For
convenience, set $\gamma(x,t)=\exp(x\zeta^{-1}+t\zeta^{-3})$, then 
\[
\psi_W=\gamma(1+a_W\zeta +O(\zeta^2)),
\]
for some function $a_W$ of $x,t$. Comparing the Fourier series for $\partial\psi_W,\psi_W$ and
$(\zeta^{-1}+\alpha^{-1})\psi_V$, all of which take values in $W$, we see that if we set
$v=a_W-a_V-\alpha^{-1}$
then 
\[
\gamma^{-1}[(\partial - v)\psi_W-(\zeta^{-1}+\alpha^{-1})\psi_V]:\{(x,t)\}\rightarrow H_+,
\]
i.e.\ it has only positive powers of $\zeta$ in its Fourier expansion. 
But $\gamma^{-1}W\cap H_+=\{0\}$ for almost all
$x,t$, hence this expression must be identically zero. Therefore
\[
(\partial -v)\psi_W=(\zeta^{-1}+\alpha^{-1})\psi_V.
\]
Further, $(\zeta^{-1}-\alpha^{-1})V=(\zeta^{-2}-\alpha^{-2})W=W$ so that a similar calculation gives
\[
(\partial + v)\psi_V = (\zeta^{-1}-\alpha^{-1})\psi_W.
\]
Therefore
\[
\begin{array}{c}
(\partial +v)(\partial - v) \psi_W=(\zeta^{-2}-\alpha^{-2})\psi_W,\\ 
(\partial -v)(\partial + v) \psi_V=(\zeta^{-2}-\alpha^{-2})\psi_V, \end{array}
\]
so in particular $u_W = -v_x-v^2+\alpha^{-2}$ and $u_V=v_x-v^2+\alpha^{-2}$. Hence if we take
$\psi=\psi_W(x,t;-\alpha)$ and $c=\alpha^{-2}$ then $(\partial^2 +u_W)\psi=c\psi$ and $(\partial
-v)\psi=0$. Thus $u_V$ is the Darboux transform of $u_W$ for the pair $(\psi,c)$.$\Box$.

\smallskip\noindent
{\bf Remark.} In fact if one examines \cite{EhlK} p5, Theorem (ii), one sees 
that a Darboux transform of a
certain type preserves the KdV spectral curve and shifts the line bundle by $\caL\mapsto\caL(Q-O)$
for some $Q\in\Sigma$ (recall that KdV solutions of finite type have spectral data
$(\Sigma,\lambda,\caL,O)$ with all the properties of the spectral data above). From this we see that 
this relationship between discrete conformal maps and Darboux transforms is true for
all periodic maps, not just those in the dressing orbit of the vacuum.

\newpage\vfill

{\large\textbf{Appendix: Non-singular spectral curves are generic.}}

\medskip\noindent
Here we will prove our earlier claim that non-singular spectral curves exist (and are therefore
generic) for periodic discrete curves with $n$ points for any $n\geq 4$.  We will use the notation 
of section 1 throughout.

Let $X_n\subset\P^1\times\ldots\P^1$ be the space of periodic discrete curves of period $n$. It is clearly
an irreducible affine open subvariety. Let $Y_n$ be the space of data $(\Sigma,O,S,P_\infty,[y])$ where:
$\Sigma$ is a complete irreducible
algebraic curve of arithmetic genus $g$ (equal to $(n-4)/2$ for $n$ even and 
$(n-3)/2$ for $n$ odd) admitting a rational function $\lambda$ of degree $2$; $O,S,P_\infty$ are
smooth points on $\Sigma$ with $\lambda$-values $0,1,\infty$ respectively, at which $\lambda$ is 
unramified (unless $n$ is odd in which case $O$ is a
ramification point); $y$ is a rational function on $\Sigma$ with divisor of poles $nS$ and $[y]$ denotes its
image in the complete linear system $\P\Gamma(\caO_\Sigma(nS))$. Since $n> 2g+2$ this linear system has
dimension $n-g$. Notice that the map $(\Sigma,O,S,P_\infty,[y])\mapsto (\Sigma,O,S,P_\infty)$ displays
$Y_n$ as a $\P^{n-g}$-bundle over the subvariety of $\P^{2g+2}$ corresponding to the possible
configurations of branch divisors. Inside $Y_n$ we consider two subvarieties: $Y_n^s$, wherein $\Sigma$
is singular; $Y_n^r$, wherein $\Sigma$ is rational with nodes only. 
In particular, $Y_n^s$ is a hypersurface in $Y$ while $Y_n^r\subset Y^s_n$ clearly
has codimension $g$ in $Y_n$.

To each $\Gamma\in X_n$ we assign the data $(\Sigma,O,S,P_\infty,[y])$ where $\Sigma,O,S,P_\infty$ are
given by the characteristic polynomial of $M_0^\lambda$ and $y=\det(H_0^\lambda)$. This $y$ has divisor
\[
D_n=\left\{\begin{array}{ll} n(O-S) & \hbox{for $n$ odd}, \\
                            \frac{n}{2}(O+\tilde O -2S)& \hbox{for $n$ even,}\end{array}\right.
\]  
and satisfies $y(P_\infty)=y(\tilde P_\infty)$. 
Thus we have an algebraic map $F:X_n\to Y_n$ with image
\[
V=\{(\Sigma,O,S,P_\infty,[y])\in Y_n:(y)=D_n, y(P_\infty)=y(\tilde P_\infty)\}.
\]
Since $V$ is irreducible, either $V\subset Y_n^s$ or
there exists a generic curve. Locally $Y^s_n$ is given by a single equation, say $Z=0$, in $Y$.
We will show that the codimension of $V\cap Y_n^r$ in $V$ is $g$ hence $Z$
cannot vanish identically on $V$ since $Y_n^r$ has only codimension $g-1$ in $Y^s_n$.

First let us describe $V\cap Y^r_n$. Each irreducible rational curve $\Sigma$ with $g$ nodes
has normalisation $\P^1$. We choose a rational parameter $t$ on $\P^1$ such that the
hyperelliptic involution is $t\mapsto 1/t$ and $t(S)=\infty$. The preimage of the singularities under the
normalisation will be $g$ pairs of the form $a_j,1/a_j$ ($j=1,\ldots g$)
and all such curves $\Sigma$ arise this way.  
For simplicity set $b=t(P_\infty)$ and $c =t(O)$. 
Then the parameters $a_j,b,c$ determine $(\Sigma,O,S,P_\infty)$. 
If $y$ has divisor $D_n$ then, up to scaling,
\[
y= \left\{\begin{array}{ll} (t-c)^n & \hbox{for $n$ odd,} \\
                            (t-c)^m(t-c^{-1})^m & \hbox{for $n=2m$,}\end{array}\right.
\]
and $a_1,\ldots,a_g,b$ must all be roots of $y(t)=y(1/t)$. Certainly this many distinct roots exist, so
$V\cap Y^r_n$ is non-empty and has dimension 1 for $n$ even, since the only free parameter is $c$, while
for $n$ odd it has dimension zero, since $O=\tilde O$ forces $c=\pm 1$. Now let us compute the 
dimension of $V$. The fibre of the map
$F:X_n\to Y_n$ over data involving any nodal rational curve is $PSL_2\times Jac(\Sigma)$ since we have
ignored M\" obius invariance and the line bundle $\caL$. Therefore $V$ has dimension $n-(g+3)$ which equals
$g+1$ for $n$ even and $g$ for $n$ odd. Therefore $V\cap Y^r_n$ has codimension $g$ in $V$
(for every $n$) whence the result follows.

\noindent
{\small Author's Addresses.
\begin{tabbing}
I Mc, P N, F P: \= GANG, Department of Mathematics \\
             \> University of Massachusetts \\
             \> Amherst, MA 01003, USA\\
             \> ian@math.umass.edu,\  norman@math.umass.edu,\\
             \> franz@gang.umass.edu \\
U H-J       \> Sfb 288, Technische Universit\" at Berlin \\
            \> Str.d.17 Juni 136\\
            \> 10623, Berlin 12, Germany\\
            \> udo@sfb288.math.tu-berlin.de\\
\end{tabbing}}

\end{document}